\documentclass[11pt,fleqn]{article}
\usepackage{hyperref}
\usepackage{amsfonts}
\usepackage{amssymb}
\usepackage{graphicx}
\newcommand{\ol}{\setlength{\itemsep}{0pt.}\begin{enumerate}}
\newcommand{\eol}{\end{enumerate}\setlength{\itemsep}{-\parsep}}
\newcommand{\ignore}[1]{}
\title{Between Arrow and Gibbard-Satterthwaite; A representation
  theoretic approach}
\author{Dvir Falik\thanks{School of Mathematical Sciences, Tel Aviv University.}, Ehud Friedgut\thanks{Einstein Institute of Mathematics, The Hebrew University of Jerusalem.
}}
\begin{document}
\date{}
\maketitle


\newtheorem{THEOREM}{Theorem}[section]
\newenvironment{theorem}{\begin{THEOREM} \hspace{-.85em} {\bf :}
}%
                        {\end{THEOREM}}
\newtheorem{LEMMA}[THEOREM]{Lemma}
\newenvironment{lemma}{\begin{LEMMA} \hspace{-.85em} {\bf :} }%
                      {\end{LEMMA}}
\newtheorem{COROLLARY}[THEOREM]{Corollary}
\newenvironment{corollary}{\begin{COROLLARY} \hspace{-.85em} {\bf
:} }%
                          {\end{COROLLARY}}
\newtheorem{PROPOSITION}[THEOREM]{Proposition}
\newenvironment{proposition}{\begin{PROPOSITION} \hspace{-.85em}
{\bf :} }%
                            {\end{PROPOSITION}}
\newtheorem{DEFINITION}[THEOREM]{Definition}
\newenvironment{definition}{\begin{DEFINITION} \hspace{-.85em} {\bf
:} \rm}%
                            {\end{DEFINITION}}
\newtheorem{EXAMPLE}[THEOREM]{Example}
\newenvironment{example}{\begin{EXAMPLE} \hspace{-.85em} {\bf :}
\rm}%
                            {\end{EXAMPLE}}
\newtheorem{CONJECTURE}[THEOREM]{Conjecture}
\newenvironment{conjecture}{\begin{CONJECTURE} \hspace{-.85em}
{\bf :} \rm}%
                            {\end{CONJECTURE}}
\newtheorem{MAINCONJECTURE}[THEOREM]{Main Conjecture}
\newenvironment{mainconjecture}{\begin{MAINCONJECTURE} \hspace{-.85em}
{\bf :} \rm}%
                            {\end{MAINCONJECTURE}}
\newtheorem{PROBLEM}[THEOREM]{Problem}
\newenvironment{problem}{\begin{PROBLEM} \hspace{-.85em} {\bf :}
\rm}%
                            {\end{PROBLEM}}
\newtheorem{QUESTION}[THEOREM]{Question}
\newenvironment{question}{\begin{QUESTION} \hspace{-.85em} {\bf :}
\rm}%
                            {\end{QUESTION}}
\newtheorem{REMARK}[THEOREM]{Remark}
\newenvironment{remark}{\begin{REMARK} \hspace{-.85em} {\bf :}
\rm}%
                            {\end{REMARK}}
\newtheorem{CLAIM}[THEOREM]{Claim}
\newenvironment{claim}{\begin{CLAIM} \hspace{-.85em} {\bf :}
\rm}%
                            {\end{CLAIM}}
\newcommand{\thm}{\begin{theorem}}
\newcommand{\lem}{\begin{lemma}}
\newcommand{\pro}{\begin{proposition}}
\newcommand{\dfn}{\begin{definition}}
\newcommand{\rem}{\begin{remark}}
\newcommand{\xam}{\begin{example}}
\newcommand{\cnj}{\begin{conjecture}}
\newcommand{\mcnj}{\begin{mainconjecture}}
\newcommand{\prb}{\begin{problem}}
\newcommand{\que}{\begin{question}}
\newcommand{\cor}{\begin{corollary}}
\newcommand{\clm}{\begin{claim}}
\newcommand{\prf}{\noindent{\bf Proof:} }
\newcommand{\ethm}{\end{theorem}}
\newcommand{\elem}{\end{lemma}}
\newcommand{\epro}{\end{proposition}}
\newcommand{\edfn}{\end{definition}}
\newcommand{\erem}{\bbox\end{remark}}
\newcommand{\exam}{\bbox\end{example}}
\newcommand{\ecnj}{\bbox\end{conjecture}}
\newcommand{\emcnj}{\bbox\end{mainconjecture}}
\newcommand{\eprb}{\bbox\end{problem}}
\newcommand{\eque}{\bbox\end{question}}
\newcommand{\ecor}{\end{corollary}}
\newcommand{\eclm}{\end{claim}}
\newcommand{\eprf}{\bbox}
\newcommand{\beqn}{\begin{equation}}
\newcommand{\eeqn}{\end{equation}}
\newcommand{\wbox}{\mbox{$\sqcap$\llap{$\sqcup$}}}
\newcommand{\bbox}{\begin{flushright} $\Box $ \end{flushright}}
\newcommand{\qed}{\bbox}
\def\sup{^}
\def\eps{\epsilon}

\def\H{\{0,1\}^n}
\def\B{\{0,1\}}

\def\S{S(n,w)}

\def\n{\lfloor \frac n2 \rfloor}
\def\ni{\noindent}
\newcommand{\rarrow}{\rightarrow}

\newcommand{\larrow}{\leftarrow}

\overfullrule=0pt
\def\setof#1{\lbrace #1 \rbrace}
\def\h#1{\widehat{ #1 }}
\def\b#1{\bar{ #1 }}
\def\hh#1{\widehat{\widehat{ #1 }}}

\def \E{\mathbb E}
\def \R{\mathbb R}
\def \Z{\mathbb Z}
\def \F{\mathbb F}
\def \C{\mathbb C}
\def \S{\mathbb S}
\def \1{\mathbf 1}

\def\<{\left<}
\def\>{\right>}
\def \({\left(}
\def \){\right)}
\def \e{\epsilon}
\def \da{\dagger}
\def \c{\cdot}
\def \l{\left}
\def \r{\right}
\def \x{\otimes}
\def \X{\bigotimes}
\def \8{\infty}

\def \z{\exp\left\{-\Omega(n)\right\}}
\def \k{O\(\exp\left\{-KN\right\}\)}
\def \A{A^{(r)}}
\def \P{{\cal P}}
\def \Q{{\cal Q}}
\def \M{{\cal M}}
\def \one{1}

\newcommand{\dnote}[1]{{ \bf (Dvir: {#1} ) }}
\newcommand{\enote}[1]{{ \bf (Ehud: {#1} ) }}
\newcommand{\dnotea}[1]{{ \bf (Dvir(1): {#1} ) }}
\newcommand{\enotea}[1]{{ \bf (Ehud(1): {#1} ) }}
\newcommand{\dnoteb}[1]{{ \bf (Dvir(2): {#1} ) }}
\newcommand{\enoteb}[1]{{ \bf (Ehud(2): {#1} ) }}

\begin{abstract}
A central theme in social choice theory is that of impossibility
theorems, such as Arrow's theorem \cite{Arr} and
the Gibbard-Satterthwaite theorem \cite{Gib,Sat}, which state that under certain natural constraints, social choice mechanisms are impossible
to construct. In recent years, beginning in Kalai \cite{Kal}, much
work has been done in finding \textit{robust} versions of these
theorems, showing ``approximate''  impossibility remains even when
most, but not all,  of the constraints
are satisfied. We study a spectrum of settings
between the case where society chooses a single outcome (\`{a}-la-Gibbard-Satterthwaite) and the choice of a complete order
(as in Arrow's theorem). We use algebraic techniques, specifically
representation theory of the symmetric group, and also prove robust
versions of the theorems that we state. Our relaxations of the
constraints involve relaxing of a version of ``independence of
irrelevant alternatives'', rather than relaxing the demand of a
transitive outcome, as is done in most other robustness results.
\end{abstract}
\section{Introduction}
Social choice deals with the aggregation of opinions of individuals in
a society into a single opinion. There are several important
\textit{impossibility} theorems in the field, stating that aggregation
mechanisms satisfying some natural conditions, are dictatorial
(dependent on the opinion of a single voter).

The first of these theorems is Arrow's theorem. Let there be a set of
$n$ individuals, who wish to decide on a ranking of $m$
alternatives. Each individual has their own full ranking of the alternatives. 
Let $L_m$ be the set of full transitive linear order on $[m]$ and
$O_m$ be the set of all anti symmetric relations on $[m]$. A social
welfare function (SWF) is a function $f:L_m^n \to O_m$, that maps the
individual rankings of the $n$ voters into an aggregated relation.
\dfn A SWF $f$ is called 
\begin{itemize}
\item \textit{Independent of Irrelevant Alternatives
(IIA)}, if for every $2$ alternatives $a,b$, the
aggregated preference between $a$ and
$b$ depends only on the individual preferences between $a$ and $b$. 
\item \textit{Consistent}, if it always returns a transitive order (is into $L_m$).
\end{itemize}
\edfn
\thm (Arrow) For $m\geq 3$, every function that is consistent and IIA,
and agrees with unanimous votes, is dictatorial. 
\ethm

Another theorem of similar flavor is Gibbard-Satterthwaite's theorem (GS), known to be strongly connected to Arrow's theorem. It deals with a setting in which the voters only wish to choose one of the $m$ alternatives. A social choice function (SCF) is a function $f:L_m^n \to [m]$, that maps the individual rankings of $n$ voters into an aggregated choice. GS deals with the game-theoretic notion of strategy proofness, where no voter has an incentive to misreport their true opinion and obtain a better result from her perspective.

For the formal definition of strategy-proofness, we introduce some notations. For a profile $x \in L_m^n, x=(x_1,...,x_n)$ and a voter $i \in [n]$, we will denote $x=(x^{-i},x_i)$, where $x^{-i}$ indicates the votes of all voters except the $i$'th. For $y \in L_m$, we will use $<_y$ to indicate the corresponding order. We similarly define $>_y, \geq_y, \leq_y$. 
\dfn
A SCF $f$ is called strategy-proof, iff 
$$\forall i \in [n], x^{-i}\in L_m^{n-1},x_i,y\in L_m, f(x^{-i},x_i) \geq_{x_i} f(x^{-i},y)$$
\edfn
i.e. no voter, under any circumstance, has an incentive to misreport
their true preference.
\thm (Gibbard-Satterthwaite) for $m\geq 3$, a social aggregator $f:L_m^n \to [m]$ that is onto and strategy-proof is dictatorial.\ethm

The connection between the notions of strategy-proofness and IIA was
demonstrated in \cite{NP,DL}, and the connection between the proofs of
these theorems was demonstrated in, e.g., \cite{Ren}. However we are
not aware of previous work which presents a single scheme that unites the different settings (SWF vs. SCF) and the different constraints (IIA vs. strategy-proofness) .

In the past decade there has been a flurry of work done in providing
analytical proofs of these theorems, and finding \textit{robust}
versions of them - i.e. showing that aggregators that almost satisfy
the constraints (consistency, IIA, strategy-proofness) are close to
fitting the classification (dictatorial). The relaxation of
consistency in Arrow's theorem was initiated in \cite{Kal} and
culminated in \cite{M}, which finally provided a robust version of the
unmodified Arrow's theorem. See also \cite{M1, Kell, Kell2}.The same
was done for several examples in the judgment aggregation setting in
\cite{Neh}. See also \cite{Xia}.

Relaxing the strategy-proofness constraint in GS has some important
computational implications. In \cite{FKKN} such a result was
achieved for functions with $m=3$, and in \cite{IKM} for
neutral functions with $m>3$. Recently  
\cite{MR} provided the final word on this theme, proving a robust version of the
unmodified Gibbard-Satterthwaite theorem. 

Our work continues this line of research. A word or two on the novelty
of our approach. The basic and beautiful idea in Kalai's paper \cite{Kal}, which
was followed in most of the subsequent work,
was using the IIA condition to translate a SWF to a set of
Boolean functions, thus enabling the application of Fourier analysis
on $\{0,1\}^n$, a technique that is very useful for proving robustness
results. The robustness referred to the distance between a given SWF
to a function which is consistent, i.e. where the output is always
transitive.  In the current paper we chose to insist on the output
being of the same form as the input (e.g. being a complete order in
the SWF setting), and measure robustness with respect to the number of violations of a special set of constraints.
 
To this end, and in order to allow a spectrum of results between SWF
(Arrow) and SCF (Gibbard-Satterthwaite), we use a notion similar to
IIA, called Independence of Rankings
(IR)  and show a robust impossibility theorem for it. IR first appeared in \cite{DH}, where they characterize impossibility domains for non-binary judgment aggregation and IR is their prime example.

We complete the analogue to Gibbard-Satterthwaite's theorem by showing that the impossibility theorem for the IR definition implies an impossibility theorem for some proper definition of strategy-proofness. This definition is an adaptation of a strategy-proofness definition introduced in \cite{DL} for binary judgment aggregation.

Our approach leads naturally to representation theory of the symmetric
group, which replaces the Abelian Fourier analysis that arises in the
previous analytical works cited above.  As in most applications of
spectral techniques to combinatorial problems, this approach includes
two components: The encoding of a combinatorial quantity as a
quadratic form, and the extraction of combinatorial information from
the spectral analysis of that form. In this work the algebraic
encoding entailed the use of block matrices. The usage of block
matrices encompasses substantial expressive power, as it enables the encoding of \textit{every} Constraint Satisfaction Problem for constraints that relate to pairs of input points, on \textit{any} size of alphabet. The canonical uses of the spectral method, utilizing standard $0-1$ matrices, are usually limited to a certain type of constraints on alphabet of size $2$, which mainly enables the treatment of notions related to expansion in graphs. In this paper the spectral analysis involves tensor algebra, allowing us to take advantage of the block structure of the matrices. This enables the extraction of the combinatorial information encoded inside the block structure.

\section{Structure of the Paper}
The paper is organized as Follows:
\begin{itemize}
\item In section \ref{sec:Results} we define the constraints we are using and state the robust impossibility theorem.
\item In section \ref{sec:Struct} we provide a bird's-eye view of the proof.
\item In section \ref{sec:Rep} we recall some essentials of representation theory.
\item In section \ref{sec:Proof} we present the formal structure of
  the proof, divided into short lemmas, which are proved in the
  subsequent section. The section is divided into two subsections,
  subsection \ref{sec:OneVoter} which deals with the case of a single
  voter, and subsection \ref{sec:NVoters}, for an arbitrary number of voters.
\item Section \ref{sec:LemProofs} provides the proofs of the lemmas from section \ref{sec:Proof}
\item The final robustness result is derived from the conclusions made in section \ref{sec:Proof} combined with an extension of a result by Friedgut, Kalai and Naor (\cite{FKN}). Section \ref{sec:FKNAdapt} provides the proof of the extension of \cite{FKN} with the necessary adaptations to our setting. The main analytic tool this proof uses is a hypercontractive inequality of Beckner and of Bonami, which is also adapted to our setting in subsection \ref{sec:Beckner}.
\item  In section \ref{sec:SP}, we introduce a strategy-proofness definition and show its connection to IR and an appropriate impossibility theorem.
\end{itemize}

\section{Results} \label{sec:Results}
In this paper we present a \textit{robust} impossibility theorem in
the flavor of Arrow's and Gibbard-Satterthwaite's theorems. The
constraint we will use is a variant of IIA. We present a single proof
dealing with functions in a spectrum of ranges, from functions
returning a full ranking (SWFs) to functions returning one alternative
(SCFs), including a plethora of ranges in between. The result can be interpreted as a $2$-query dictatorship test with full completeness.

In our setting, we deal with aggregation of rankings of $m$ alternatives. A ranking is a permutation $x\in\S_m$. We will use the convention $x(rank)=name$.

For presentation sake, we shall begin with the definition of the constraint when used for functions returning a full ranking (SWFs) and state the corresponding impossibility theorem without robustness. The more complicated definitions and theorems will follow.

\subsection{Main Theorem}
\dfn \label{dfn:IRManyVtrs} A social aggregator $f:\S_m^n \to \S_m$ satisfies Independence of Rankings (IR) iff the aggregated ranking of the $j$'th alternative is dependent only on the individual rankings of the $j$'th alternative
$$\forall x,y\in\S_m^n, j \in [m], \(\forall i \in [n], x_i^{-1}(j)=y_i^{-1}(j)\)\Rightarrow f(x)^{-1}(j)=f(y)^{-1}(j)$$
\edfn

This constraint requires independence of rankings instead of independence of pairwise preferences required in IIA. This constraint was discussed in \cite{DH}, in the context of non-binary judgment aggregation.

 As in IIA, this definition compares voting profiles which may differ
 in any number of votes. Throughout this paper, we shall use an
 alternative, equivalent definition, that compares inputs that differ
 in a single vote (as is the case in the definition of strategy-proofness):

\dfn \label{dfn:IR1Vtr} A social aggregator $f:\S_m^n \to \S_m$ satisfies Independence of Rankings (IR) iff
$$\forall i \in [n], j \in [m], x^{-i} \in\S_m^{n-1}, x_i,y_i\in\S_m, x_i^{-1}(j)=y_i^{-1}(j)\Rightarrow f(x^{-i},x_i)^{-1}(j)=f(x^{-i},y_i)^{-1}(j)$$
\edfn
It is easy to show that these two definitions are equivalent, via a hybrid argument. The corresponding impossibility theorem is
\thm \label{thm:MainSWF} For $m\geq 3$, a social aggregator $f:\S_m^n \to \S_m$ that is IR is either a constant function or dictatorial of the following form: there exists a voter $i$ and a constant permutation $y$ of the rankings such that $f(x)=y\circ x_i$.\ethm

In \cite{DH}, a similar impossibility theorem was shown, using purely
combinatorial arguments. Their result deals with a larger range of
possible formats of input and output, yet demands further constraints
on the function in question. We do not see a way to extend their
techniques to achieve robustness.

\subsection{Robust Impossibility Theorem}
A robust impossibility theorem means that when the constraint is
\textit{almost} satisfied, then the function is \textit{almost}
dictatorial. 

To the best of our knowledge, most previous robustness results
regarding SWF's focused on relaxation of the rationality (the
transitivity of the outcome), and measure the distance to a function
which is rational. We, instead, demand rationality, and relax (our
variant) of IIA.  There is
something satisfying about relaxing this specific constraint, rather
than others, as it seems to be slightly less natural than rationality and unanimity.

\dfn A social aggregator $f:\S_m^n \to \S_m$ is called $\epsilon-IR$
if the rate of constraints that are  \textbf{not} satisfied is smaller
than $\eps$, i.e.
$$\sum_{i \in [n],j \in [m]} Pr_{ x^{-i} \in\S_m^{n-1},
  x_i,y_i\in\S_m}  \left[\left( x_i^{-1}(j)=y_i^{-1}(j)\right) \wedge\
  \left( f(x^{-i},x_i)^{-1}(j)\neq f(x^{-i},y_i)^{-1}(j)\right)\right] \leq \epsilon$$
\edfn

\thm \label{thm:MainSWFRbst} For $m\geq 3$, a social aggregator
$f:\S_m^n \to \S_m$ that is $\epsilon$-IR is $O(m^8\epsilon))$
close to a function that is either a constant function or dictatorial
of the following form: there exists a voter $i$ and a constant
permutation $y$ of the rankings such that $f(x)=y\circ x_i$. \ethm

\subsection{A Spectrum of Ranges}
As stated earlier, we will also deal with a setting where the
aggregated opinion is not a full ranking, but a partial ranking. Let
$H\subseteq \S_m$ be a subgroup of $\S_m$. We call it a
\textit{fixing} subgroup if it consists of all permutations respecting
a given partition of the $m$ rankings into $2$ or more parts. An $H$-social
aggregator is a function $f:\S_m^n\to\S_m/H$, where $\S_m/H$ refers to
right cosets of $H$ (We use this notation even though $H$ is not a
normal subgroup). Many types of functions fall under this scheme. Examples are:
\begin{itemize}
\item For $H$ as the trivial group, $H$-social aggregators are SWFs.
\item For $H$ as the group of permutations fixing the element $1$, $H$-social aggregators are SCFs.
\item For $H$ as the group of permutations fixing the set $\{1,2,3\}$, $H$-social aggregators are functions returning triumvirates.
\item For $H$ as the group of permutations fixing the sets $\{1\}$ and $\{2,3\}$, $H$-social aggregators are functions returning a president and two vice-presidents.
\end{itemize}
The definition of IR can be extended to $H$-social aggregators in the following manner:
\dfn Let $H\subseteq \S_m$ be  a fixing subgroup of $\S_m$. For $H_1$
a right coset of $H$ in $\S_m$, and $j \in [m]$ define the $j$-profile of $H_1$ as the multiset $H_1^{-1}(j)=\left\{y^{-1}(j) | y \in H_1\right\}$.
\edfn
\dfn \label{dfn:H-IR} Let $H\subseteq \S_m$ be a fixing subgroup of $\S_m$. An $H$-social aggregator $f$ satisfies Independence of Rankings (IR) iff the aggregated $j$-profile is dependent only on the individual rankings of the $j$'th alternative.
$$\forall i \in [n], j \in [m], x^{-i} \in\S_m^{n-1}, x_i,y_i\in\S_m,
x_i^{-1}(j)=y_i^{-1}(j)\Rightarrow
f(x^{-i},x_i)^{-1}(j)=f(x^{-i},y_i)^{-1}(j)$$
\edfn
\xam
Let $H=\S_{m,1}$, the group of permutations that fix the element $1$. In that case, an $H$ social aggregator $f$ is a function that returns a single winner in an election.

When the winner of the election is $k$, then the function returns the coset $H^k$, which is the coset that includes all permutations that assign $k$ to $1$. The $k$-profile of $H^k$ is $|H|$ copies of $1$. The $j$-profile of $H^k$ for every $j \neq k$ has $\frac{|H|}{m-1}$ copies of each number between $2$ and $m$.

$f$ is IR in that case if, for an alternative $j$, when given the individual rankings of $j$ by all voters in a voting profile $x$, we are able to determine whether $j$ is the winner of the election (i.e., the $j$-profile of $f(x)$ is $|H|$ copies of $1$) or not (i.e. the $j$-profile of $f(x)$ has $\frac{|H|}{m-1}$ copies of each number between $2$ and $m$).

\exam
We shall leave the exact definition of an $\epsilon$-IR $H$-social
aggregator to a later part of the paper, see definition \ref{dfn:H-IRf}.  The impossibility theorems also extend to $H$-social aggregators.

\thm \label{thm:MainGen} Let $H\subseteq \S_m$ be a fixing subgroup of
$\S_m$. For $m\geq 3$, an $H$-social aggregator $f$ that is IR is
either a constant function or dictatorial of the following form: there
exists a voter $i$ and a constant permutation $y$ of the rankings such
that $f(x)=Hy\circ x_i$.\ethm

\thm \label{thm:MainGenRbst} Let $H\subseteq \S_m$ be a fixing subgroup of $\S_m$. For $m\geq 3$, an $H$-social aggregator $f$ that is $\epsilon$-IR is $O_H(poly(m)\epsilon)$ close to a function that is either a constant function or dictatorial of the following form: there exists a voter $i$ and a constant permutation $y$ of the rankings such that $f(x)=y\circ x_i$.\ethm

\section{Structure of the proof}\label{sec:Struct}
We give here a short exposition of the proof. For simplicity, we shall refer here to the basic form of the Main theorem (theorem \ref{thm:MainSWF}), where the function is a SWF. To simplify the notation, we shall also use in this section definition \ref{dfn:IRManyVtrs} for IR, instead of \ref{dfn:IR1Vtr}, which is the definition we shall use in the rest of the paper.

We shall treat this problem as a constraint satisfaction problem (CSP). We shall use the definition \\
\texttt{
\\\indent\indent Find all functions $f:\S_m^n \to \S_m$ s.t. \\
\indent IR: $j \in [m], x,y \in \S_m^n, $\\
\indent\indent $x^{-1}(j)=y^{-1}(j) \Rightarrow (f(x))^{-1}(j)=(f(y))^{-1}(j)$\\}

\ni A CSP has a generic algebraic encoding. The function $f$ can be encoded as a function returning a vector in $\R^{\S_m}$, which is the characteristic vector of the singleton $\{f(x)\}$. This encoding can be interpreted as a tensor $F \in \R^{\S_m^n\times\S_m}$, with $2$ indices $x,v\in\S_m$  $$F_{x,v}=\1_{v=f(x)}.$$

\ni The constraints can be algebraically encoded using a matrix that
represents their truth table, or, in our case, since we want to count
the number of violated constraints, the truth table of their negation. We use a matrix of matrices. For every two inputs $x,y \in \S_m^n$, the $(x,y)$'th entry of the matrix will be a matrix in $\R^{\S_m \times \S_m}$. This matrix will be the truth table of the negation of the constraints concerning $x,y$ and $j$. This encoding can also be interpreted as a tensor: 
$$\(\(L^j\)_{xy}\)_{v_xv_y}=\1_{[x^{-1}(j)=y^{-1}(j)]
  \wedge [v_x^{-1}(j)\neq v_y^{-1}(j)]}$$
 where $j \in [m], x,y \in \S_m^n$ and $v_x,v_y\in\S_m$.

\ni We can use these tensors in a quadratic form to count the number
if unsatisfied constraints. Since $L$ is the truth table of the
negation of the constraints, the quadratic form $$\sum_j FL^jF^t$$
counts the number of 
violated constraints.

The CSP under this encoding takes the form\\
\texttt{
\\\indent\indent Find all $F \in \R^{\S_m^n\times\S_m}$ s.t. \\
\indent Consistency: $\forall x \in \S_m^n, F_{f{x}*}$ is a characteristic vector of a singleton\\
\indent IR: $\sum_j FL^jF = 0$\\}

\ni The proof unfolds as follows:
\begin{itemize}
\item We show that $L \succeq 0$ (PSD), i.e. $\sum_j FL^jF^t \geq 0$ for every $F$. Therefore, the functions that satisfy IR are precisely the kernel of L.
\item Explicitly find the kernel of $L$, using diagonalization.
\item Show that all consistent functions in the kernel of $L$ are dictatorships.
\end{itemize}

\ni As for the robustness, we will show that functions that are $\epsilon$-IR are $L_2$ close to the kernel of $L$. We shall generalize the result of \cite{FKN} to prove that such functions, that are also consistent, are $L_2$ close to dictatorships. 

\ni For an $H$ social aggregator, we shall encode $f$ to return characteristic vectors of cosets of $H$, normalized so that their $L_1$ norm equals $1$. We shall call such vectors $H$ coset vectors. As a tensor $F$, this encoding takes the form:
$$F_{x,v}=\frac{1}{|H|}\1_{v\in f(x)}$$ 
A very convenient feature of our definitions and approach is that the
introduction of $H$ social aggregators does not insert any new elements to the
proof. Essentially,
the same Laplacian $L$ encodes the notions of IR and $\epsilon$-IR for
$H$ social aggregators, and only the consistency constraint changes. The algebraic CSP for $H$ social aggregators is
\texttt{
\\\indent\indent Find all $F \in \R^{\S_m^n\times\S_m}$ s.t. \\
\indent Consistency: $\forall x \in \S_m^n, F_{\mathbf{x}*}$ is an $H$-coset vector\\
\indent IR: $\sum_j FL^jF^t = 0$\\}

\section{Representation Theory} \label{sec:Rep}
In this section we recall some basic notions of representation theory
that are necessary for our proof.

A representation is a Homomorphism $\rho$ from a group $G$ to $GL_d(\C)$, the group of complex d-dimensional square matrices. $d$ is called the dimension of the representation $d(\rho)$. A representation is called irreducible if it is not similar to a direct sum of $2$ representations.

For a finite group, there is a one-to-one correspondence between the
conjugacy classes of the group and irreducible representations (up to
similarity). We shall denote the number of conjugacy classes of $G$
as $[G]$. For a conjugacy class $k \in [[G]]$, its corresponding
irreducible representation will be denoted as $\rho^k$. We will
sometimes consider $\rho^k$ as a function, and sometimes treat it as a
vector with $|G|$ entries all of which are  $d(\rho^k)\times d(\rho^k)$ matrices.

In this paper we will deal exclusively with the symmetric group or direct products of the symmetric group.
It is well known that in the symmetric group one can choose a basis
for which all irreducible representations (also known as irreps) have real, unitary matrices
as values. Henceforth we will assume we are dealing with such a basis.

The defining representation of the symmetric group $\S_m$ is the
permutation representation $P$ of dimension $m$. 
$$P(x)_{ij}=\1_{x(i)=j}$$
It is well known that $P=\rho^0\oplus\rho^1$, i.e. it is the direct sum of two irreducible
representations: the trivial one, which we denote by $\rho^0$, and the
$(n-1)$-dimensional $\rho^1$. Specifically, if one chooses a basis for
which $\rho^0, \rho^1$ are real and unitary then there exists an
orthonormal $m\times m$ matrix $U$ such
that 
\begin{equation}\label{eqn:PDec}
P(x)=U\(\rho^0(x)\oplus\rho^1(x)\)U^t
\end{equation}
where here the $\oplus$ refers to a matrix composed of blocks.

The all ones vector spans the one dimensional eigenspace of $\R^m$ corresponding to
the trivial representation component of $P$, hence $U$ can be written
in the following form (where $C$ is a $m\times(m-1)$ matrix, and the $C_i$'s are its rows)
\begin{equation} \label{eq:C}
U=\(\begin{array}{ll}{\begin{array}{c}\frac 1{\sqrt{m}}\\ \vdots \\ \frac 1{\sqrt{m}}\\ \end{array}}& C\end{array}\)=\(\begin{array}{ll}\frac 1{\sqrt{m}} & C_1\\ \vdots & \vdots \\ \frac 1{\sqrt{m}} & C_m\\ \end{array}\)
\end{equation}

The character of a representation $\rho$, denoted by $\chi_\rho$, is the trace of the representation: $\chi_\rho(x)=tr(\rho(x))$. It is easy to see the the character of similar representations $\rho$ and $U\rho U^{-1}$ are the same.

A final tool we wish to recall is Schur's orthogonality, which states that the vectors of the form $\rho^{\bf k}_{ij}$ are orthogonal, (but not necessarily orthonormal) 
$$\sum_x \rho^{ k_1}(x)_{i_1j_1}\rho^{ k_2}(x)_{i_2j_2} = \delta_{ k_1k_2}\delta_{ i_1i_2}\delta_{ j_1j_2} \frac {m!}{d(\rho^{? k_1})}$$ For a finite group $G$, these vectors form a complete orthogonal basis for the set of functions from $G$ to $\C$.

Schur's orthogonality implies that the characters of irreducible representations are orthonormal: $<\chi_{\rho^k},\chi_{\rho^l}>=\delta_{k,l}$. This implies that for a representation $\tau$ and an irreducible representation $\rho$, $<\chi_{\tau},\chi_{\rho}>$ is the multiplicity of $\rho$ in the decomposition of $\tau$ to irreducible representations. 

A complete set of irreducible representations for $G^n$ is the set of tensors of the irreps of $G$:
$$\{\rho^{\b r}=\X_{i=1}^n \rho^{r_i}\}_{\b r \in [[G]]^n}$$

\subsection{Fourier transform and diagonalization} \label{sec:Fou}
For a finite group $G$, given a function $f:G \to \C$, its Fourier transform at a representation $\rho$ is 
$$\h{f}(\rho)=\E_{x \in G} f(x)\rho(x).$$
 We shall sometimes use the abbreviated notation $\h{f}(k)=\h{f}(\rho^k)$

For such a function $f$, define a matrix $M\in\C^{G \times G}$
whose values are $M_{x,y}=f(x^{-1}y)$, then $M$ can be
partially diagonalized (decomposed into eigenspaces) using representation theory. In a partial diagonalization, an eigenspace does not have a corresponding eigenvalue, but rather a corresponding \textit{eigenblock}. Given a chosen set of basis vectors $\{v_i\}_i$ for an eigenspace of a matrix A, the $(i,j)$'th entry of the corresponding eigenblock is $v_iAv_j^t$.

In the case of $M$, each irrep $\rho^k$ defines an eigenspace spanned by the aforementioned orthogonal vectors $\(\rho_{i,j}(x)\)_{x \in G}$.
The eigenblock corresponding to $\rho^k$ is $I_{d(\rho^k)} \x \h{f}(k)$. 

Another formulation of this is: $M$ is a $|G| \times |G|$ matrix,
whose entries are given by 
$$M_{x,y}=\sum_{k \in [[G]]}d(\rho^k) tr\(\(\b{\rho^k}\h{f}(k){\b{\rho^k}}^t\)_{x,y}\)$$
Where $\b{\rho^k}$ is a column  vector of size $|G|$ whose $x$'th
entry is the matrix $\rho^k(x)$. The multiplication
$\h{f}(k){\b{\rho^k}}^t $ means multiplying the {\em entries} of  $
{\b{\rho^k}}^t$ by $\h{f}(k)$, hence $\b{\rho^k}\h{f}(k){\b{\rho^k}}^t$ is a $|G|$ dimensional matrix whose entries are $d(\rho^k)$ dimensional matrices.

For a function $g$, $gM$ is known as the convolution of $g$ and $f$ and is denoted as $g*f$. The partial diagonalization discussed above shows that:
\begin{eqnarray}
\nonumber\h{g*f}(k)&=&\h{g}(k)\h{f}(k)\\
\nonumber (gMg^t)_{x,y}&=&\sum_{k \in [[G]]}d(\rho^k)tr\(\(\h{g}(k)\h{f}(k)\h{g}^t(k)\)_{x,y}\)
\end{eqnarray}

If $f$ is a characteristic function of a set $T$ of generators of $G$, then $M$ is the adjacency matrix of the Cayley graph $\Gamma(G,T)$. For convenience, we shall call such a graph a Cayley graph even when $f$ is a characteristic function of any subset $T$ of $G$, as the property of $T$ generating $G$ is irrelevant for our uses.

\section{The Proof} \label{sec:Proof}
In this section we present a more detailed version of the proof, divided into lemmas. The actual proofs of the lemmas will appear in  section \ref{sec:LemProofs}.
 
\subsection{One Voter Functions}\label{sec:OneVoter}
We begin our analysis by treating the case of a single voter, since
this contains the analytical seed from which the multi-voter case
grows. The combinatorial problem for a single voter is not very
difficult, although if one is interested in a robustness theorem it
seems that straightforward elementary techniques are insufficient.

In order to treat a social welfare functions on one voter
$f:\S_m\to\S_m$, we construct in this section a quadratic form
encoding the constraints, and, as mentioned in section
\ref{sec:Struct}, we will diagonalize it.  

For $f:\S_m \to \S_m$, denote by $IR(f)$ the rate of unsatisfied
constraints $$IR(f)=\sum_{j \in [m]} 
Pr_{x,y \in \S_m}  [ \left(x^{-1}(j)=y^{-1}(j)\right) \wedge \left(f(x)^{-1}(j)\neq f(y)^{-1}(j)\right) ].$$

 Let $X^j$ be the matrix $X\in\R^{\S_m\times\S_m}$, 
$$X^j_{xy}=\1_{x^{-1}(j)=y^{-1}(j)}=\1_{x^{-1}y(j)=j}.$$
Let $\bar{X^j}$ be its complement $\bar{X^j}_{xy}=1-X^j_{xy}$.
We will use the vector encoding described in \ref{sec:Struct}, $F_{x,v}=\1_{v=f(x)}$ (or the corresponding definition for $H$-social aggregators).

The following lemma describes a quadratic form in the values of $f$ that equals $IR(f)$:
\lem \label{lem:LapStr}
Let $f$ be a social aggregator and $F$ its encoding as described above. Let $L'^j= X^j\x \bar{X^j}$, and $L'=\sum_j L'^j$, then
$$IR(f)= \frac{1}{|\S_m|^2} FL'F^t$$
\elem

$X^j$ is the adjacency matrix of the Cayley graph 
$\Gamma(\S_m,\S_{m,j})$, where $\S_{m,j}$ is the subgroup of $S_m$ of permutations fixing $j$. A quadratic form based on the Laplacian of that same graph is more suitable for our purposes, because it is PSD. The Laplacian of that graph, $Y$, is given by 
$$Y^j=(m-1)!I-X^j$$
The corresponding quadratic form is given in the following lemma. The quadratic forms given in lemmas \ref{lem:LapStr} and \ref{lem:LapLap} are equivalent when $F$ represents a consistent function.
\lem \label{lem:LapLap}
Let $f$ be a social aggregator and $F$ its encoding as described above. Let $L''^j= Y^j \x X^j$, , and $L''=\sum_j L''^j$, then
$$IR(f)= \frac{1}{|\S_m|^2} FL''F^t$$
\elem

Since $X$ the adjacency matrix of a Cayley graph, it can be partially diagonalized
via the representations of the symmetric group, as explained in subsection \ref{sec:Fou}. As will be shown in the proof of the following lemma, all the information relevant to the computation of $IR(f)$ lies in $X$'s $\rho^1$ component. This leads to a simplified quadratic form, used with a different encoding for $f$. Let $g:\S_m^1 \to \R^{(m-1)\times(m-1)}$ be a an encoding of $f$ such that $g(x)=\rho^1(f(x))$. A vector form of $g$ is a vector $G$ whose each entry is a $m-1\times m-1$ matrix $G_{x}=g(x)$. (For $H$ social aggregators, $g(x)=\E_{y\in f(x)}\rho^1(y)$). 

The corresponding quadratic form is  as follows.
\lem \label{lem:LapFou}
Let $f$ be a social aggregator and $G$ its encoding as described above. Let $L^j= Y^j \x D^j$, where $D^j=C_j^tC_j$ (See \ref{eq:C}), and $L=\sum_j L^j$, then
\begin{equation} \label{eqn:LapFou}
  IR(f)=\frac{1}{|\S_m|^2} tr (GLG^t)  
\end{equation}
\elem

For an $H$ social aggregator, $IR(f)$ is defined as such: it averages, for pairs of inputs that agree on the ranking of $j$, the square of the $\ell_2$ distance of the characteristic vector of the $j$-profile of the outputs.
\dfn \label{dfn:H-IRf} For an $H$-social aggregator $f:\S_m^n\to\S_m/H$, define $IR(f)$ to be
$$IR(f)=\sum_{i,j}\E_{x^{-i},x_i,y_i} \(\1_{x_i^{-1}(j)=y_i^{-1}(j)}\left\|\frac{{\bf n}_{f(x^{-i},x_i)^{-1}(j)}-{\bf n}_{f(x^{-i},y_i)^{-1}(j)}}{|H|}\right\|_2^2\)$$
Where for a multiset $S$, ${\bf n}_S$ is its characteristic vector, i.e., for an element $x$, ${\bf n}_S$ at the index $x$ equals the number of occurrences of $x$ in $S$.

\edfn

In the following we claim that the quadratic form we have defined in \ref{lem:LapFou} indeed equals the value $IR(f)$, as defined in definition \ref{dfn:H-IRf}. As usual, we only address the $1$ voter case in this subsection:
\clm \label{clm:LapFou+H-IR}
Let $H\subseteq \S_m$ be a fixing subgroup of $\S_m$. Let $f$ be an  $H$-social aggregator over $1$ voter, and $G$ be its encoding as defined above, then 
$$IR(f)=\frac{1}{|\S_m|^2}tr(GLG^t)$$
\eclm

We partially diagonalize $L$ in the following lemma, decomposing it to
eigenspaces. For the purposes of this lemma, we shall define the operator $\tilde{tr}$ that operates on block matrices. The operator returns a matrix whose $(x,y)$'th entry is the trace of $(x,y)$'th block in the original matrix: $\tilde{tr}(M)_{xy}=tr(M_{xy})$. 

\lem \label {lem:LDiag}
$$L=\sum_{r\in [[\S_m]]} d(\rho^r) \tilde{tr}\(\( \rho^r\ \x I\) \h{L}(r) \(\rho^r \x I\)^t\)$$
where
\begin{eqnarray}
\nonumber \h{L}(0)&=&I\cdot 0\ \ ,\ \ \h{L}(r>1)=I\x I \cdot \frac 1m\\
\nonumber \h{L}(1)&=&\frac 1{m-1} \(\frac{m-1}{m} I\x I-  \sum_j D^j \x D^j\)
\end{eqnarray}
\elem

The diagonalization of $\h{L}(1)$ is given by:
\lem\label{lem:hatL1Diag}
$\h{L}(1)$ has $3$ orthogonal eigenspaces whose dimensions are $1,m-1,(m-1)^2-m$.  Denote their corresponding basis matrices as $U_0,U_1,U_2$. The corresponding eigenvalues are $0,\frac{1}{m(m-1)},\frac{1}{m}$. The eigenvectors are vectors in $\R^{(m-1)\c (m-1)}$. When read as a $(m-1)\times(m-1)$ matrix, $U_0$ is the identity matrix.
\elem

This diagonalization proves that $L$ is PSD and hence, in light of
Lemma \ref{lem:LapFou}, determines that its kernel is the space of IR functions. This is summarized in this corollary:
\cor \label{cor:1VtrIRfnc}
Let $f:  \S_m \to \S_m/H$ be an IR function,
and let $g:\S_m \to \R^{m-1\times m-1}$ be defined by $g(x)= \rho^1(f(x))$ then there exist vectors $a$ and $b$ in $\R^{(m-1)\cdot(m-1)}$, that can be
read as $(m-1)\times(m-1)$ matrices $A$ and $B$) such that 
$$g_{x}=b\rho^0(x) + \tilde{tr}\((a U_0^t)(I\x \rho^1(x)\) = B + A\rho^1(x)$$
\ecor
To complete the characterization of $1$-voter IR functions, we present the following claim and lemma. We shall not prove the lemma as it is a consequence of the $n$ voter case.
\clm \label{clm:FixSbgrp} Let $H$ be a subgroup of $\S_m$, and let $f$ be an $H$-social aggregator and $g$ be $g(x)=\rho^1(f(x))$, then $$\forall x, g(x)g^t(x)=M$$ When $H$ is a fixing subgroup, then $M \neq 0$.
\eclm
\lem \label{lem:1VtrCons}
Let $g:\S_m \to \R^{m-1\times m-1}$ be a function of the form $g(x)=B+ A\rho^1(x)$ that satisfies the constraint $\forall x, g(x)g^t(x)=M$ for some constant matrix $M\neq 0$. Then either $B=0$ or $A=0$.
\elem

\subsection{Many Voter Functions}\label{sec:NVoters}
The quadratic form capturing the notion of $IR(f)$ for functions on
$n$ voters is constructed using the quadratic form for $1$ voter, in
the following lemma.
 
\lem\label{lem:LapNVtr}
Let $H$ be a fixing subgroup of $\S_m$. For a function
$f:\S_m^n\to\S_m/H$, let $G$ be as before, the encoding of $f$ in
terms of its $\rho^1$ component:
$$G(x)=\rho^1(f(x))=\E_{y\in f(x)}\rho^1(y),$$
where the expectation with respect to $y$ refers to the case of
$H$-social-aggregators. 
Let 
$$L^{n,j,i}=  I^{\x i-1} \x Y^j \x I^{\x n-i}\x D^j\ \ ,\ \ L^{n,i}=\sum_j L^{n,j,i}\ \ L^n=\sum_i L^{n,i}$$
Then the number of unsatisfied constraints is 
$$IR(f)=\frac{1}{|\S_m|^{n+1}} tr(GL^{n}G^t)$$
\elem
We can diagonalize $L^n$ based on our diagonalization of $L$:
\cor \label{cor:LapNVtrDiag} The diagonalization of $L^n$ is given by:
$$L^{n,i}=\frac 1{|\S_m|^{n-1}}\sum_{\bar{r}\in[[\S_m]]^n} d(\rho^{\bar{r}})\tilde{tr} \left( (\rho^{\bar r} \x I )\h{L^{n,i}}({\bar{r}}) (\rho^{\bar{r}}\x I)^t  \right)$$
The $\h{L^{n,i}}$'s are derived from the 1 voter $\h{L}$'s, as follows.
$$
\(\h{L^{n,i}}(\rho^{\b r})\)_{k_1...k_nkl_1...l_nl}=\prod_{t \in [n],t\neq i} \(I_{d\(\rho^{r_t}\)}\)_{k_tl_t} \c \(\h{L}(\rho^{r_i})\)_{k_ikl_il}
$$


\ecor

The $\h{L^n}$ coefficients are matrices which are not necessarily diagonal. The following lemma partly characterizes their diagonalization, in a manner that suffices for our needs. 

\lem \label{lem:LapNVtrDiag2}
\begin{enumerate}
\item If there exists any coordinate $i \in [n]$ for which $r_i>1$ then $\h{L^n}(\b r) \succeq \frac 1m I$ ($A \succeq B$ means that $A-B$ is PSD).
\item Otherwise, if there exist at least $2$ coordinates $i \in [n]$ for which $r_i=1$ then $\h{L^n}(\b r) \succeq O\(\frac 1 {m^2}\) I$
\item Otherwise, if there exists exactly $1$ coordinate $i \in [n]$ for which $r_i=1$ then $\h{L^n}(\b r)=\h{L}(1)$. As shown in lemma \ref{lem:hatL1Diag}, it has a $0$ eigenvalue corresponding to the eigenvector $U_0$ and a smallest nonzero eigenvalue of $O\(\frac 1 {m^2}\)$.
\item Otherwise, $\b r$ is all zeros, and $\h{L^n}(\b r)=0$.

\end{enumerate}
\elem

For a PSD matrix, its \textit{spectral gap} is its smallest non-zero eigenvalue. We conclude the following from the diagonalization of $L^n$:
\cor \label{cor:NVtrIRfnc}
\begin {itemize}
\item The kernel of $L^n$, which is the set of all IR functions, consists exclusively of functions of the form $g(x_1,...,x_n)=B+\sum_{i=1}^nA^i\cdot\rho^1(x_i)$
\item The spectral gap of $\frac 1{|\S_m|}L^n$ is $\frac 1{O(m^2)}$. 
\end{itemize}
\ecor

To finish the proof of theorem \ref{thm:MainSWF}, we need to show that the intersection of the kernel of $L^n$ with the consistency constraint, includes only dictatorships. We don't need to use the consistency constraint to its full capacity. All we need to use is the quadratic constraint that $\forall x, g(x)g^t(x)=M$, where $M$ is some constant matrix. In the proofs section, we shall show that this constraint is valid for any $H$. For instance, if $f$ is a SWF ($H$ is the trivial group), then since $\rho^1$ is unitary, $\forall x, g(x)g^t(x)=I$. 

\cor \label{cor:IR+cons}
IR functions which are consistent are dictatorships.
\ecor

\subsection{Robustness}
By using the quadratic form, we are able to connect the combinatorial notion of $IR(f)$ to the analytical notion of the distance between $f$ and the kernel of $L^n$. This connection depends on the spectral gap of $L$:
\cor \label{cor:RobSpecGap}
If $IR(f)\leq\eps$, and $g$ is the encoding of $f$ as above, then there exists a function $h$ in the kernel of $L^n$ such that $\|h-g\|_2^2 \equiv \E_{x\in \S_m^n} \|h(x)-g(x)\|_2^2 \leq O(m^2)\eps$
\ecor

In corollary \ref{cor:IR+cons} we characterized $IR$ functions, using the fact that the function satisfies two constraints:
\begin{itemize} 
\item Being IR, which a linear constraint, because it is equivalent to being in the kernel of $L$.
\item Being consistent, which a quadratic constraint.
\end{itemize}
It is clear that the intersection of a linear and a quadratic constraint can contain only a few points.
Indeed, we showed that consistent functions which satisfy the linear constraint are dictatorial.

For $\eps$-IR functions, corollary \ref{cor:RobSpecGap} shows that the first (linear) constraint is relaxed to being $L_2$ close to the linear constraint. We wish to show that consistent functions that are $L_2$ close to the linear constraint are $L_2$ close to being dictatorial.

In \cite{FKN} a similar result was shown for Boolean functions on Boolean variables. It was shown that Boolean functions that are linear are dictatorial, and that Boolean functions that are $L_2$ close to being linear are $L_2$ close to a dictatorial function. Being Boolean is, naturally, a quadratic constraint. We adapt this theorem to  our setting. From it we deduce our main theorem.

\section{Proofs for section \ref{sec:Proof}} \label{sec:LemProofs}
\subsection{Proofs of the Lemmas for 1 voter functions}
\ni\textbf{Proof of lemma \ref{lem:LapStr}:} This is the straightforward definition of the anti-constraints.

$$\frac{1}{|\S_m|^2}\sum_{jxv_xyv_y}F_{xv_x}L'_{xv_xyv_y}F_{yv_y}=\frac{1}{|\S_m|^2}\sum_{jxv_xyv_y} F_{xv_x}X^j_{xy}\bar{X^j}_{v_xv_y}F_{yv_y}=$$$$\frac{1}{|\S_m|^2}\sum_{jxv_xyv_y} \1_{v_x=f(x)}\1_{x^{-1}(j)=y^{-1}(j)}\1_{v_x^{-1}(j)\neq v_y^{-1}(j)}\1_{v_x=f(x)}=$$$$\sum_{j}\E_{xy} \1_{x^{-1}(j)=y^{-1}(j)}\1_{f(x)^{-1}(j)\neq f(y)^{-1}(j)}$$\qed

\ni\textbf{Proof of Lemma \ref{lem:LapLap}:} 

Recall 
$$
\begin{array}{lclcccclccc}
L'' & = & \sum_j & X^j &\x &\(J-X^j\)& = &\sum_j &X^j \x J &-&X^j\x X^j\\
L' &= &\sum_j &\((m-1)!I-X^j\)&\x &X^j=&\sum_j &(m-1)!I\x X^j&-&X^j \x X^j\\
\end{array}
$$
Therefore, we need to show that
\begin{equation} \label{eqn:L'vsL''}
F\(\sum_j X^j\c J\)F^t = F\(\sum_j(m-1)!I\x X^j\)F^t
\end{equation}
when $F$ is consistent. 

Since $F$ is consistent, $\forall x, \sum_{v}F_{xv}=1$, so the left hand side of (\ref{eqn:L'vsL''}) is 
$$\sum_{jxyv_xv_y} F_{xv_x}\(X^j_{xy}\c 1\)F_{yv_y}=\sum_{jxy} X^j_{xy}=m\c m! \c (m-1)!=m!^2$$
Since the diagonal of $X^j$, for every $j$, is all ones, the right hand side of (\ref{eqn:L'vsL''}) is
$$\sum_{jxyv_xv_y} F_{xv_x}\((m-1)!\delta_{xy}\x X^j_{v_xv_y}\)F_{yv_y}=\sum_{jxv_xv_y} F_{xv_x}\((m-1)! X^j_{v_xv_y}\)F_{xv_y}=$$$$\sum_{jx} \((m-1)! X^j_{f(x)f(x)}\)=\sum_{jx} \((m-1)!1\)=m\c m!\c(m-1)!=m!^2$$

\ni\qed

Before we carry on, this is good place to recall $U$ and $C$ from
equation \ref{eq:C}. The orthonormality of $U$ implies the following:
\clm
\begin{eqnarray} \label{eq:CTen}
\nonumber CC^t&=&I-\frac Jm\\
\nonumber C^tC&=&I\\
\nonumber \1 C&=&0
\end{eqnarray}
\eclm

\ni \textbf{Proof of Lemma \ref{lem:LapFou}:} 
We will need two simple claims and their corollary.

\ni{\bf Claim:} $$X^j_{xy}=\(U\(\rho^0(xy^{-1})\oplus\rho^1(xy^{-1})\)U^t\)_{jj}$$

\ni{\bf Proof:} Recall $P$, the defining representation of $\S_m$. A permutation $x$ has a fixed point $j$ iff $(P_x)_{jj}$ is $1$. Therefore,
$$X^j_{xy}=(P_{x^{-1}y})_{jj}$$

Recall that $P_{x}=U(\rho^0\oplus\rho^1)U^t$, and that $U$ is orthonormal. Therefore
$$X^j_{xy}=(P_{xy^{-1}})_{jj}=\(U\(\rho^0(xy^{-1})\oplus\rho^1(xy^{-1})\)U^t\)_{jj}$$
\qed

\ni{\bf Claim:} 
\begin{equation} \label{eqn:XDiag}
    X^j_{xy}=
 \(\frac{\1 \rho^0(xy^{-1})\1^t}m + C \rho^1(xy^{-1})C^t\)_{jj} 
\end{equation}
\ni{\bf Proof:} Follows from the expansion of $U$.\qed
Next denote $X^j=X^{j,0}+X^{j,1}$, where
$X^{j,0}_{xy}=
 \(\frac{\1 \rho^0(xy^{-1})\1^t}m\)_{jj} $ and
$X^{j,1}_{xy}=
 \(C \rho^1(xy^{-1})C^t\)_{jj} $.
Likewise, denote $L^j=L^{j,0}+L{j,1}$ where $L^{j,0}=Y^j\x X^{j,0}$ and $L^{j,1}=Y^j\x x^{j,1}$.

\ni{\bf Corollary:}
$$IR(f)=FL'F^t=\sum_j \(FL^{j,0}F^t + FL^{j,1}F^t\) = \sum_j FL^{j,1}F^t$$

\ni{\bf Proof:}
The first equalities follow from the expansion of $IR(f)$ and $L'$. We now show that $FL^{j,0}F^t=0$. 

$$FL^{j,0}F^t=\sum_{xyv_xv_y}\frac{1}{|S_m|^2}F_{xv_x}\( Y^j_{xy} \c\( \frac {\1\rho^0(v_x)\rho^0(v^{-1}_y)\1^t}m\)_{jj}\)F_{yv_y}$$
Since $F$ is consistent, and $\rho^0_x=\one_x$, we have $\sum_{v_x}F_{xv_x}\rho^0(v_x)=\1_x$. Therefore
$$FL^{j,0}F^t=\sum_{xy}\frac{1}{|S_m|^2}\1_x\(Y^j_{xy}  \frac {1_j1_j}m\)\1_y=0$$ \qed

\ni It follows that 
$$IR(f)=\sum_j FL^{j,1}F^t=\sum_{jxyv_xv_y}\frac{1}{|S_m|^2}F_{xv_x}\(Y^j_{xy} \( C\rho^1(v_x)\rho^1(v^{-1}_y)C^t\)_{jj}
\)F_{yv_y}$$
In light of the above we can define $G$, a vector whose entries are $(m-1)$ dimensional matrices, as $G_{x}=\sum_{v}F_{xv}\rho^1(v)$, and the quadratic form becomes
$$IR(f)=\sum_{jxykl_xl_y}\frac{1}{|S_m|^2}\(G_{x}\)_{kl_x}\(Y^j_{xy}  C_{jl_x}C_{jl_y}
\)\(G_{y}\)_{kl_y}=\sum_{j}\frac{1}{|S_m|^2}tr\(G\(Y^j \x D^j\)G^t\)$$

This completes the Proof of Lemma \ref{lem:LapFou}.
\qed

\ni\textbf{Proof of claim \ref{clm:LapFou+H-IR}:}
Recall definition \ref{dfn:H-IRf}:
$$IR(f)=\sum_{j}\E_{x,y} \(\1_{x^{-1}(j)=y^{-1}(j)}\left\|\frac{{\bf n}_{f(x)^{-1}(j)}-{\bf n}_{f(y)^{-1}(j)}}{|H|}\right\|_2^2\)$$

We will show the equivalency via two simple claims:
\begin{itemize}
\item[{\bf Claim a:}] $$\frac{1}{|\S_m|^2}tr(GLG^t)=\sum_{j}\E_{x,y} \(\1_{x^{-1}(j)=y^{-1}(j)}\|C^jg^t(x)-C^jg^t(y)\|_2^2\)$$
\item[{\bf Claim b:}] $$\left\|\frac{{\bf n}_{f(x)^{-1}(j)}-{\bf n}_{f(y)^{-1}(j)}}{|H|}\right\|_2^2=\|C^jg^t(x)-C^jg^t(y)\|_2^2$$
\end{itemize}
\ni{\bf Proof of Claim a:} 
Quadratic forms based on matrices which are Laplacians of graphs, such as $Y^j$, naturally differentiate values across edges of the graph. This can be seen via the following simple decomposition of $Y^j$: Define the matrix $Z^{j,x,y}$ to be the Laplacian of the graph whose vertices are the elements of $\S_m$ and is either the empty graph or contains a single edge connecting $x$ with $y$ in case $x^{-1}(j)=y^{-1}(j)$.

Clearly, $Y^j=\sum_{(x,y)\in{\S_m \choose 2}}Z^{j,x,y}$. It is also easy to see that the diagonalization of the $Z$'s is given by $Z^{j,x,y}=\1_{xy^-1(j)=j} {d^{x,y}}^td^{x,y}$ where $d^{x,y}$ is a vector that has $1$ in $x$, $-1$ in $y$ and $0$ otherwise.

Using this, we get
$$tr\(GL^jG^t\)=\sum_{(x,y)\in{\S_m \choose 2}}tr\(G\(Z^{j,x,y}\x D^j\)G^t\)=$$$$\sum_{(x,y)\in{\S_m \choose 2}}\1_{xy^-1(j)=j} tr\(G\({d^{x,y}}^t\x{C^j}^t\)\c\({d^{x,y}}\x{C^j}\)G^t\)=$$$$\sum_{(x,y)\in{\S_m \choose 2}}\1_{xy^-1(j)=j} tr\( \(g(x)-g(y)\){C^j}^t\c{C^j}\(g(x)-g(y)\)^t \)=$$$$\sum_{(x,y)\in{\S_m \choose 2}}\1_{xy^-1(j)=j} \< {C^j}\(g(x)-g(y)\)^t,{C^j}\(g(x)-g(y)\)^t \>$$ \qed

\ni{\bf Proof of claim b:}
Clearly, the normalized characteristic vector of the $j$-profile of $g(x)$ is $\frac 1{|H|}{\bf n}_{f(x)^{-1}(j)}=e_j\(\E_{y\in g(x)}P^t_y\)$, where $e_j$ is the $j$'th unit vector. When transforming this vector using the orthonormal matrix $U$, we get:
$$\frac 1{|H|}{\bf n}_{f(x)^{-1}(j)}U=\(e_j\E_{y\in g(x)}P^t_y\)U = \(e_j U\(1\oplus g^t(x)\)U^t\)U=e_j U\(1 \oplus g^t(x)\)=$$$$\(\frac 1{\sqrt{m}}\ \ C_j\)\(1\oplus g^t(x)\)=\(\frac 1{\sqrt{m}}\ \ C^jg^t(x)\)$$
Therefore, since $U$ is orthonormal, 
$$\left\|\frac{{\bf n}_{f(x)^{-1}(j)}-{\bf n}_{f(y)^{-1}(j)}}{|H|}\right\|_2^2=\left\|\frac{{\bf n}_{f(x)^{-1}(j)}-{\bf n}_{f(y)^{-1}(j)}}{|H|}U\right\|_2^2=$$$$
\left\|\(\frac 1{\sqrt{m}}\ \ C^jg^t(x)\)-\(\frac 1{\sqrt{m}}\ \ C^jg^t(y)\)\right\|_2^2=\|C^jg^t(x)-C^jg^t(y)\|_2^2$$

\qed

\ni\textbf{Proof of lemma \ref{lem:LDiag}:} 
Recall that $$L=\sum_j \((m-1)!I - X^j\) \x D^j$$
The identity matrix can be decomposed using the irreps of $\S_m$ and
Schur orthogonality.
$$(m-1)!I_{xy}=\frac{(m-1)!}{m!}\(\sum_{r \in [[\S_m]]}
d(\rho^r)\rho^r(x)\rho^r(y^{-1})\)$$
Recall the decomposition of $X^j$ from (\ref{eqn:XDiag}): 
$$X^j_{xy}=
\frac {\rho^0(x)\rho^0(y^{-1})}m+ C_{j}\rho^1(x)\rho^1(y^{-1})C^t_{j}
$$

Summing these decompositions, and using the fact that $\sum_j D^j=\sum_j C^t_{j}C_{j}=C^tC=I$ yields that for a fixed pair
  $x,y$, one has

$$\begin{array}{rccrccclc}
L_{xy}&=&&\rho^0(x)\x I&\c& \frac 1m I - \frac 1m I &\c &   \rho^0(y^{-1})\x I&+\\
      &&(m-1)\tilde{tr}(&\rho^1(x)\x I &\c& \frac 1{m-1} \(\frac{m-1}{m} I\x I-  \sum_j D^j \x D^j\) &\c &\rho^1(y^{-1})\x I&)+\\
      & &\sum_{k>1}^{|[\S_m]|}d(\rho^k)\tilde{tr}(&\rho^k(x)\x I &\c& \frac 1m I\x \(\sum_jC^t_{j}C_{j}\) &\c &\rho^k(y^{-1})\x I&)\\
\hspace*{0pt}\\
&=&&\rho^0(x)\x I&\c& 0 &\c &   \rho^0(y^{-1})\x I&+\\
      &&(m-1)\tilde{tr}(&\rho^1(x)\x I &\c& \frac 1{m-1} \(\frac{m-1}{m} I\x I-  \sum_j D^j \x D^j\) &\c &\rho^1(y^{-1})\x I&)+\\
      & &\sum_{k>1}^{|[\S_m]|}d(\rho^k)\tilde{tr}(&\rho^k(x)\x I &\c& \frac 1m I\x I &\c &\rho^k(y^{-1})\x I&)
\end{array}
$$
as required.
\qed

\textbf{Proof of lemma \ref{lem:hatL1Diag}:}
Denote by $E$, a $m\times (m-1)^2$ matrix whose $j$'s row is $C_j\x C_j$. We need to diagonalize $$\sum_j D^j\x D^j=\sum_j \(C^t_jC_j\)\x\(C^t_jC_j\)=\sum_j\(C^t_j\x C^t_j\)\(C_j\x C_j\)=E^tE$$

\ni The nonzero eigenvalues of $E^tE$ are the nonzero eigenvalues of
$EE^t$ (this can be deduced from the SVD decomposition of E).
Recall $CC^t=I-\frac{J}m$. Therefore, $$\(EE^t\)_{ij}=\(C_i\x C_i\)\(C^t_j\x C^t_j\)=\(C_iC^t_j\)^2=\(CC^t\)^2_{ij}=$$$$\(\delta_{ ij}-\frac 1m\)^2=\(1-\frac 2m\)\delta_{ij}+\(\frac 1m\)^2$$
Therefore, $EE^t=\frac {m-2}m I+\frac J{m^2}$, and its eigenvalues are $\frac{m-1}m$ with multiplicity $1$ and $\frac{m-2}m$ with multiplicity $m-1$. 

We can verify that the eigenvector of $E^tE$ corresponding to the $\frac{m-1}m$ eigenvalue is $U_0$, which is the identity matrix parsed as a vector. We need to use two simple facts: 
\begin{itemize}
\item
For $3$ matrices $A,B$ and $C$, the term $A\c B\c C^t$ when parsed as a vector, is equal to the term $\(A\x C\)\c B$, when $B$ is parsed as a vector.
\item
For a matrix $C$ whose rows are $\{C_j\}_j$, the term $\sum_j C_j \x C_j$ equals to the term $C^tC$ parsed as a row vector. This is because $C^tC=\sum_j C^t_jC_j$.
\end{itemize}
Therefore, we get that
$$(E\c U_0)_j=(C_j\x C_j)U_0=C_jIC^t_j=C_jC^t_j=\(CC^t\)_{jj}=\(I-\frac Jm\)_{jj}=\frac {m-1}m$$
Which means that $EU_0=\1 \frac{m-1}m$ and

$$\(E^tE\)U_0=E^t\1 \c \frac{m-1}m=\frac{m-1}m \(\sum_j C_j \x C_j\)^t$$
Since $\sum_j C_j \x C_j$ is $C^tC$ parsed as a vector and  $C^tC=I$ and $U_0$ is $I$ parsed as a vector, we get that
$\(E^tE\)U_0=\frac{m-1}m U_0$.


\qed

\ni\textbf{Proof of corollary \ref{cor:1VtrIRfnc}:} We have shown that eigenspaces of $L$ with eigenvalue $0$ are $\rho^0$ and  $\rho^1U^0$. All the other eigenvalues are positive, so $L$ is PSD, and IR functions are in the kernel of $L$.

Therefore, if $G$ is an IR function:
\begin{enumerate}
\item Its $\rho^0$ Fourier coefficient can be anything. 
\item Its  $\rho^1$ Fourier coefficient must be of the form $A\(U^0\)^t$ for some vector $A$. Recall that this Fourier coefficient has $4$ indices, and that $U^0$ is the identity matrix, parsed as a vector, and also $A$ is some matrix parsed as a vector. 
\item All of its other Fourier coefficients must be $0$
\end{enumerate}
Explicitly, $$g(x)_{kl}=B_{kl}\rho^0(x) + \sum_{ts}\rho^1_{ts}(x)\delta_{sl}A_{kt} = B_{kl} + \sum_t A_{kt}\rho^1_{tl}(x)$$
This means that the function $g$ is $g(x)=B+A\cdot \rho^1(x)$ (where
$A$ and $B$ are parsed as $(m-1)\times(m-1)$ matrices), so $g$ is a linear
function in $\rho^1(x)$, and, according to lemma \ref{lem:1VtrCons}, since
$g$ is consistent,
 either $A$ or $B$ are $0$.


\qed

{\bf Proof of claim \ref{clm:FixSbgrp}:} Let $M_H$ be $M_H=\E_{x\in H} \rho^1(x)$. Since $g$ is consistent, it must be of the form $\forall x, \exists y, g(x)=M_H\rho^1(y)$. Therefore, 
$$ g(x)g(x)^t=M_H\rho^1(y)\rho^1(y)^tM_H^t=M_HM_H^t=\E_{x\in H} \rho^1(x) \E_{y\in H} \rho^1(y^{-1})=$$$$\E_{x,y\in H} \rho^1(xy^{-1})=\E_{z\in H} \rho^1(z)=M_H$$.

We now turn to show that $M_H \neq 0$ when $H$ is fixing. Denote $P_H=\E_{x\in H} P(x)$. Clearly $P_H=U(1\oplus M_H)U^t$. Therefore,  if $M_H=0$, then $P_H \simeq J$. However, let $i$ and $j$ be $2$ indices from $2$ different parts of the partition that $H$ fixes, then clearly $(P_H)_{ii}\neq(P_H)_{ij}$, because for every $x\in H$, $P_{ij}(x)=0$, and for some $x\in H$, $P_{ii}(x)=1$. 
\qed

\rem Notice that when $H$ is not a fixing subgroup, it might be the case that $M_H=0$, and then $g \equiv 0$, and every consistent $H$ social aggregator satisfies IR. This is the case, for instance, when $H$ is the group of even permutations.
\erem

\subsection{Proofs of Lemmas for many voter functions}
\ni\textbf{Proof of lemma \ref{lem:LapNVtr}:}
We show that applying the $n$ voter quadratic form is similar to applying the $1$ voter quadratic form but over pairs of inputs that differ in only one voter: 
$$GL^{n,i}G^t=\sum_{x,y}G(x)\(I^{\x i-1}\x Y^j \x I^{\x n-i}\x D^j \) G^t(y)=$$$$
\sum_{x,y}G(x)\(\prod_{t=1}^{i-1}I_{x_t,y_t} \x Y^j \prod_{t=i+1}^{n}I_{x_t,y_t}\x D^j \) G(y^{-1})=$$
$$
\sum_{x^{-i},x_i,y_i}G(x^{-i},x_i,y_i)\( Y^j \x D^j \) G^t(x^{-i},y_i)$$
From here, the proof follows exactly in the path of the proof of claim \ref{clm:LapFou+H-IR}.
\qed

\ni\textbf{Proof of lemma \ref{lem:LapNVtrDiag2}:}
\begin{enumerate}
\item
Let $i$ be such that $r_i>1$, then, by corollary \ref{cor:LapNVtrDiag}
$$\h{L^{n,i}}(\b r)=\frac 1m \X_{k=1}^n I_{d(r_k)} \x I_{m-1}$$
Since $\h{L^{n}}(\b r)=\sum_j \h{L^{n,j}}(\b r)$ and for every $j, \b r$, $\h{L^{n,j}}(\b r) \succeq 0$, we have $$\h{L^{n}}(\b r) = \h{L^{n,i}}(\b r) + \sum_{j\neq i} \h{L^{n,j}}(\b r) \succeq \h{L^{n,i}}(\b r)$$
\item We shall focus only on the case where there are exactly $2$
  distinct $i$ and $j$ such that $r_i=r_j=1$, as this case
  produces the minimal eigenvalue.  Indeed, assume that there are $k>2$ such indices $i_1,i_2,...i_k$ for which $r_{i_j}=1$, and denote $r'$ to be $1$ in $i_1$ and $i_2$ and $0$ otherwise, then $\h{L^{n,i_1}}(r)+\h{L^{n,i_2}}(r)=\h{L^n}(r')\x I_{m-1}$. Since the rest of the terms $\h{L^{n,i_k}}(r)$ are PSD, The minimal eigenvalue of $\h{L^n}$ is at least as large as the minimal eigenvalue of $\h{L^n}(r')$.



Recall the diagonalization of $\h L (1)$ (from lemma \ref{lem:hatL1Diag}) $$\h{L}(1)\succeq \frac 1{m(m-1)} \(I\x I- \frac 1{m-1} U^0{U^{0}}^t\)$$ From this, it is easy to deduce the following:
\begin{equation} \label{eqn:2rho1eigval}
    \h{L^n}(\b r)=\h{L^{n,i}}(\b r)+\h{L^{n,j}}(\b r)\succeq \frac 1{m(m-1)}  \(2I\x I\x I- \frac 1{m-1}\(A^i{A^{i}}^t+A^j{A^{j}}^t\)\)
\end{equation}
Where $A^i$ is a $(m-1)^3\times (m-1)$ matrix of the form
$A_{(pqs)t}=\delta_{pq}\delta_{st}$ ($p,q,s$ and $t$ are indices going
from $1$ to $m-1$. $(pqs)$ forms the row index of $A^i$ and $t$ is its
column index). Likewise, $A^j$ is a $(m-1)^3\times (m-1)$ matrix of
the form $A_{(pqs)t}=\delta_{ps}\delta_{qt}$. Denote the right hand
side of (\ref{eqn:2rho1eigval}) as Q.

Denote $B=A^i+A^j$ and $C=A^i-A^j$. It is easy to see that
$A^i{A^i}^t+A^j{A^j}^t = \frac 12\(BB^t+CC^t\)$.
The following are easy to verify:
\begin{itemize}
\item $B^tC=C^tB=0$
\item $BB^t=2m\c I$
\item $CC^t=(2m-4)\c I$

\end{itemize}
Therefore, we may deduce that the columns of $B$ and the columns of $C$ are orthogonal eigenvectors of $A^i{A^i}^t+A^j{A^j}^t$ with eigenvalues $m$ and $m-2$, respectively. Plugging this into the expression for $Q$, we get that the minimal eigenvalue of $Q$, corresponding to the columns of $B$, is $\frac 1{m(m-1)} (2-\frac{m}{m-1})=\frac {m-2}{m(m-1)^2}$.

\ni\textbf{Important Note:} Notice that if $m=2$, $Q$ has eigenvalues equal to $0$, and therefore we cannot deduce that the function is a dictatorship for $m=2$.

\ni Items {\bf 3} and {\bf 4} are trivial.

\end{enumerate}

\textbf{Proof of corollary \ref{cor:IR+cons}:}
We need to show that only one of $B,A_1,...,A_n$ is not zero. We begin by showing that w.l.o.g., we may assume that $\E_x g(x)=0$ and
  therefore $B=0$. Indeed, we introduce a dummy variable $y \in \S_m$
  and define
 $$g'(x,y)=g(xy^{-1})\rho^1(y)$$ 
Where for $x\in \S_m^n$ we denote $xy = (x_1y,x_2y,...,x_ny)$. Note that  $\E_{x,y}g'(x,y)= 0$ because $$\E_{x,y}g'(x,y)=\E_{x,y}g'(xy,y)=\E_{x,y}g(x)\rho^1(y)=\E_xg(x)\E_y\rho^1(y)=0$$ and that $IR(g')=IR(g)$. Assume the claim is true for $g$ such that $\E g=0$, apply it to $g'$ to get that either $g'(x,y)=A_i\rho^1(x_i)$ for some $i$, or $g'(x,y)=B\rho^1(y)$. In the first case, it follows that $g(x)=A_i\rho^1(x_i)$, and in the second case, it follows that $g(x)=B$. 

We now assume that $B=0$. Recall that by claim \ref{clm:FixSbgrp}, we have $g(x)g(x)^t=M$ for some $M\neq 0$. On the other hand,
$$g(x)g(x)^t=\sum_{ij}A^i\rho^1(x_i){\rho^{1}}^t(x_j){A^j}^t$$


The summand for $i,j$, translates to 
$$\(A^i\rho^1(x_i){\rho^{1}}^t(x_j){A^j}^t\)_{tu}=\sum_{p,q,w,s}\rho^1_{pq}(x_i)\rho^1_{ws}(x_j)\(A^i_{tp}A^j_{wu}\delta_{qs}\)$$

From this expansion we may deduce the Fourier coefficient of $gg^t$ at $\b r$ where $\b r$ is $1$ at $i$ and $j$ and $0$ otherwise:
\begin{equation} \label{eqn:fHatDict}
  \h{gg^t}_{tpwuqs}(\b r)=(m-1)^2\(\(A^i_{tp}A^j_{wu}\delta_{qs}\)+\(A^j_{tp}A^i_{wu}\delta_{qs}\)\)  
\end{equation}

Since $gg^t$ is a constant function and the Fourier expansion is unique, we get that $\h{gg^t}(\b r)$ must be $0$.
Assume by contradiction that there exist indices $tp$ and $wu$ where
$A^i_{tp}\neq 0,A^j_{wu}\neq 0$. Equation (\ref{eqn:fHatDict}) implies that $A^j_{tp}A^i_{wu}=-A^i_{tp}A^j_{wu}\neq 0$. Therefore, $A^j_{tp}\neq 0$ and $\h{gg^t}_{tptpqq}(\b r)\simeq A^i_{tp}A^j_{tp}\neq 0$

\qed

\textbf{Proof of corollary \ref{cor:RobSpecGap}: }

The result follows from the following general statement regarding quadratic forms.
\\ \ni{\bf Claim:} Let $K$ be a finite set. Equip the linear space $\R^K$ with the $L_2$ metric $d(u,v)=\E_{x \in K}(v_x-u_x)^2$. Let $M\in \R^{K \times K}$ be a PSD matrix with spectral gap $\lambda$, then for any vector $u$,
$$\frac 1{|K|} uMu^t \leq \lambda d(u,ker(M))$$
where $ker(M)$ is the kernel of $M$ and $d(u,ker(M))$ is the minimal distance of $u$ from an element of $ker(M)$ according to $d$.

\ni{\bf Proof:} Assume that $M$ has eigenvalues $0=\lambda_1\leq \lambda_2\leq...\leq\lambda_{|K|}$, with corresponding eigenvectors $v_1,...v_{|K|}$, and that $s$ is the first index for which $\lambda_s>0$. Denote by $\h{u}_i$ the projection of $u$ on $v_i$. Clearly, the  distance of $u$ from the kernel of $M$ is $\frac 1{|K|} \sum_s^{|K|} \h{u}_i^2$, and the value of the quadratic form on $u$ is $$u^tMu=\sum_1^{|K|} \h{u}_i^2 \lambda_i = \sum_s^{|K|} \h{u}_i^2 \lambda_i \leq \lambda_s \(\sum_s^{|K|} \h{u}_i^2\)$$ \qed

In our case, we have shown in lemma \ref{lem:LapNVtr} that $IR(f)=\frac 1{|\S_m|^{n+1}}tr(GLG^t)$. When the
elements of $G$ are parsed as vectors instead of $(m-1) \times (m-1)$
matrices, this translates to 
$$IR(f)=\frac 1{|\S_m|^{n+1}} G\cdot \(I\x L\) G^t$$
Since the spectral gap of $\frac 1{|\S_m|}L$ is $O\(\frac1 {m^2}\)$, so is the spectral gap of $\frac 1{|\S_m|}I\x L$.

Therefore, if $IR(f)\leq \epsilon$, then there exists $h$ in the kernel of $I\x L$ such that $\|G-h\|_2\leq O(m^2)\epsilon$ (when the elements of $G$ and $h$ are parsed as vectors).
\qed

\section{Adapted version of \cite{FKN}} \label{sec:FKNAdapt}

For the sake of self-containedness, we present here the proof from FKN with minor modifications, needed for the application of the theorem to our setting.



\ni The adapted theorem goes as follows
\thm 
\label{thm:dict}
Let $Lin(\S_m^n)$ be the space of functions of the form $\sum_i A^i\rho^1(x_i)$. Let $g:\S_m^n \to \R^{m-1 \times m-1}$ be a function such that 
\begin{itemize}
\item $\E g =0$
\item There exists a matrix $M$ such that  $tr(M)=1$ and $\forall x\in \S_m^n,
  g(x)g^t(x)=M$.
\item $\|g-Lin(\S_m^n)\|_2^2\leq \tau$.
\end{itemize}
then there exists $i$ such that $\E\|g-A^i\rho^1(x_i)\|_2^2 \leq O(m^5\tau)$.   
\ethm

Before we carry on with the proof of theorem \ref{thm:dict}, we shall show how to apply it to get the proof of the main theorem:

\textbf{Proof of theorem \ref{thm:MainGenRbst}:}  Let $f$ be an $\eps$-IR H social aggregator and let $g$ be $g(x)=\rho^1(f(x))$, as before. We assume, w.l.o.g., that $\E g=0$ (see the proof of \ref{cor:IR+cons}). Since $g$ is consistent, by claim \ref{clm:FixSbgrp} it holds that there exists some matrix $M$ such that for every x, $g(x)g^t(x)=M$. Let $K=tr(M)$. By corollary \ref{cor:NVtrIRfnc} we have $\E \|g-Lin(\S_m^n)\|_2^2\leq O(m^2)\eps$. Therefore, we may apply theorem \ref{thm:dict} on $\frac 1{\sqrt K} g$ and get that there exists $i$ such that $\E\|g(x)-A^i\rho^1(x_i)\|_2^2 \leq O(Km^7\eps)$. This is almost what we need, except that we have no guarantee that $A^i\rho^1(x_i)$ is consistent. However, we show in lemma \ref{lem:Round} that we can round it to a consistent function without losing more than a multiplicative constant factor in the distance. 

The value of $K$ ranges between $1$ (for SCF's) and $m$ (for SWF's).\qed  

Next we show, as promised, that we can "round off" our approximation to an approximation which is consistent without
losing too much. 
\lem \label{lem:Round}Let $g:\S_m^n \to \R^{(m-1)\times(m-1)}$ be a function such that there exists a matrix $M$ such that $range(g)\subseteq M\c range(\rho^1)$. Assume that there exists a function $h:\S_m^n \to \R^{(m-1)\times(m-1)}$, where $h(x)=A\rho^1(x_i)$ for some $i$, and $\|g-h\|_2\leq \delta$, then there exists a function $h'=A'\rho^1(x_i)$ such that $range(h')\subseteq M\c range(\rho^1)$, and $\|g-h'\|_2\leq 2\delta$
\elem
\prf
Assume that there does not exist $A'\in M\c range(\rho^1)$ such that $|A-A'|_2\leq\delta$, then $$|g-h|_2=|g\rho^1(x_i^{-1})-h\rho^1(x_i^{-1})|_2=|g\rho^1(x_i^{-1})-A|_2\geq d(M\c range(\rho^1),A)>\delta$$
which is a contradiction. Therefore, there exists such $A'$. Define $h'$ as above $h'(x)=A'\rho^1(x_i)$, then
$$|g-h'|_2\leq|g-h|_2+|h-h'|_2\leq 2\delta$$
\eprf 

\textbf{Proof of theorem \ref{thm:dict}:}

Let $h$ be the projection of $g$ on $Lin(\S_m^n)$ and let $$q=g-h$$
and
$$r= f\c f^t-M$$ 
Since $\E g=0$, there exist $A^i$'s such that $h=\sum_i A^i\rho^1(x_i)$.
We will show that since $h \c h^t$ is close to $M$, $r$ is typically
close to $0$, and we will deduce some information on the $A^i$'s.

In the proof of corollary \ref{cor:IR+cons}, we used the fact that when $\eps=0$, $r\equiv 0$.
For the case when $\eps$ is positive, we will try and show that $r$ is
close to the $0$ function.

\rem \label{rem:rQuad}
Note that the entries of $r$ are of the form
$$ r_{k,l}(x)=\(C +\sum_{i,j}A_i\rho^1(x_i){\rho^1}^t(x_j)A_j^t\)_{k,l}$$
Which makes functions of degree $2$. For a formal definition of the degree of a function, see definition \ref{dfn:deg}, to come in subsection \ref{sec:Beckner}
\erem


\rem In the following lemma, we use a constant $C=\frac
1m^{\frac 12}$ that is defined in corollary \ref{cor:BeckFinal}. We
shall call it $C$ from now on so that it will be easier to trace its
origin.
\erem
\begin{lemma}
\label{Rsmall}
$$E(r^2) \leq K\eps=108(m-1)^4C^8\eps.$$
\end{lemma}
\cor
There exists   $i$ such that
$\|A^i\|_2^2 \geq 1 -
\left(1+ \frac {(m-1)K} {1 - \eps}\right)\eps=1-O(m^3)\eps.$
\ecor
Clearly, this corollary implies theorem \ref{thm:dict}.

{\bf Proof of corollary:}

Because of the orthogonality of $h$ and $q$, and because $\E \|g\|_2^2=1$, $\E \|h\|_2^2 = 1-\eps$.

Also, because of the orthogonality of the $\rho^1(x_i)$'s,
$$\E \|h\|_2^2=\sum_i\E\|A^i\rho^1(x_i)\|_2^2=\sum_i\|A^i\|_2^2$$

$$\E \|r\|_2^2 \geq \sum_{i\neq j} \E\|A^i\rho^1(x_i)(A^j\rho^1(x_j))^t\|_2^2$$

Expanding this expression we get: 
$$\E\|A^i\rho^1(x_i)(A^j\rho^1(x_j))^t\|_2^2=$$$$\E_{x_1x_2}\sum_{\alpha\mu\beta_1\gamma_1\eta_1\beta_2\gamma_2\eta_2}A^i_{\alpha\beta_1}\rho^1_{\beta_1\gamma_1}(x_i)\rho^1_{\eta_1\gamma_1}(x_j)A^j_{\mu\eta_1}A^i_{\alpha\beta_2}\rho^1_{\beta_2\gamma_2}(x_i)\rho^1_{\eta_2\gamma_2}(x_j)A^j_{\mu\eta_2}=$$$$
\frac{1}{(m-1)^2}\sum_{\alpha\mu\beta_1\gamma_1\eta_1\beta_2\gamma_2\eta_2}A^i_{\alpha\beta_1}\delta_{\beta_1\beta_2}\delta_{\gamma_1\gamma_2}\delta_{\eta_1\eta_2}\delta_{\gamma_1\gamma_2}A^j_{\mu\eta_1}A^i_{\alpha\beta_2}A^j_{\mu\eta_2}=$$$$
\frac{1}{m-1}\sum_{\alpha\mu\beta\eta}A^i_{\alpha\beta}A^j_{\mu\eta}A^i_{\alpha\beta}A^j_{\mu\eta}=$$$$
\frac{1}{m-1}\|A^i\|_2^2\|A^j\|_2^2$$

Therefore, if $\max \|A^i\|_2^2 = t$
$$\(\E \|h\|_2^2\)^2 = (1-\eps)^2=\sum_i\(\|A^i\|_2^2\)^2 + \sum_{i,j,i\neq j}\|A^i\|_2^2\|A^j\|_2^2\leq $$$$(1-\eps)t + (m-1)\E \|r\|_2^2 \leq (1-\eps)t + (m-1)K\eps$$

which gives the desired bound on $t$.

{\begin{flushright} $\Box $ \end{flushright}}

{\bf Proof of Lemma \ref{Rsmall} : } The proof
consists of
 two parts: first we will show that typically $r$ obtains values
 close to 0. Then we will use a hypercontractive estimate due
 to Beckner and Bonami to bound higher moments of $r$ in terms of
 its second moment showing that its tail decays fast enough.
\begin{lemma}
\label{p}
 Let $0< \alpha < 1/4$ be a constant to be chosen later.
Let
$$p=Prob(|r|_2^2> \alpha^2).$$
Then
$$p \leq \frac{16 (m-1)^2\eps}{\alpha^2}.$$
\end{lemma}
We defer the proof of this lemma for the moment.
\begin{lemma}
\label{less}
$$E(\|r\|_2^2) \leq \frac{(m-1)^2\alpha^2}{1-4(m-1)C^4\sqrt{\eps}/\alpha}.$$
\end{lemma}

Choosing the optimal value of $\alpha$ (which is $6(m-1)C^4\sqrt{\eps}$)
immediately proves Lemma \ref{Rsmall}.
So, to finish the proof we
now present the proofs of Lemmas \ref{p} and \ref{less}.
\\ {\bf Proof of Lemma \ref{p}:}
Recall that $q=g-h$ and that $\E\|q\|_2^2=\eps$. Using the fact that  $gg^t=M$ yields
$$r = hh^t-M = qq^t-gq^t-qg^t$$
By the triangle inequality,
$$\|r(x)\|_2 \leq \|r(x)\|_1 \leq \|q(x)q^t(x)\|_1+2\|g(x)q^t(x)\|_1$$
We shall prove in claim \ref{clm:MatCS} that, for $d \times d$ matrices $A$ and $B$, $\|AB\|_1\leq d\|A\|_2\|B\|_2$. This implies  
$$\|r(x)\|_2 \leq (m-1)(\|q(x)\|^2_2+2\|g(x)\|_2\|q(x)\|_2)$$

denote $t=\|q(x)\|_2$. We use that fact that $\|g(x)\|_2^2=tr\(g(x)g^t(x)\)=tr(M)=1$ 
$$ \|r(x)\|_2 \leq (m-1)(2t +t^2) $$
A simple analysis shows that if $t < \alpha/4(m-1)$,
then $\|r(x)\|_2 < \alpha$. Hence by Markov's inequality
$$Pr[\|r(x)\|_2 > \alpha] \leq Pr[ \|q(x)\|_2^2 > (\alpha/4(m-1))^2] \leq
\frac{(4(m-1))^2 \eps}{\alpha^2}.$$
{\begin{flushright} $\Box $ \end{flushright}}
\clm \label{clm:MatCS}
The following is the version of Cauchy-Schwartz that we used in the
previous proof. Let $A$ and $B$ be two $d$ dimensional real matrices, then
$$\|AB\|_1\leq d\|A\|_2\|B\|_2$$
\eclm
\prf 
Let $J$ be the all ones matrix. It is easy to see that for a matrix $A$, $\|A\|_1=<A,J>$. $J$ can be decomposed to a sum of $d$ permutation matrices $P^1,...P^d$. Then,
$$ \|AB\|_1=<AB,J>=\sum_{i=1}^d<AB,P^i>=\sum_{i=1}^d<A,B^tP^i>\leq\sum_{i=1}^d
\|A\|_2\|B^tP^i\|_2=d\|A\|_2\|B\|_2$$
When $A=B=J$, this inequality is tight.
\eprf
 {\bf Proof of Lemma \ref{less}:} For convenience of notation let
$X=E\|r\|_2^2$, $X_{ij}=E(r_{ij}^2)$ and $Y_{ij}=E(r_{ij}^4)$. Also denote $p_{ij}=P(|r_{ij}|\leq \alpha)$. In corollary \ref{cor:BeckFinal} we show a Beckner type inequality  of the form: $$Y_{ij} \leq C^8 X_{ij}^2$$ (Notice that, as mentioned in remark \ref{rem:rQuad}, the entries of $r$ are real functions of degree $2$ with Fourier coefficients only in $\rho^0$ and $\rho^1$, and therefore conform to the requirements of corollary \ref{cor:BeckFinal})
Using this we obtain

$$ X = \sum_{ij} X_{ij}= \sum_{ij}E(r_{ij}^2) = \sum_{ij}\((1-p_{ij})E(r_{ij}^2|r_{ij}^2 \leq \alpha^2) + p_{ij}E(r_{ij}^2| r_{ij}^2 >
\alpha^2)\) \leq$$
$$\sum_{ij}\( (1-p_{ij})\alpha^2 + p_{ij} \sqrt{E(r_{ij}^4| r_{ij}^2 > \alpha^2)} \)\leq$$
$$\sum_{ij} \(\alpha^2 + p_{ij} \sqrt{\frac{Y_{ij}}{p_{ij}}} \)\leq $$
$$\sum_{ij} \(\alpha^2 + \sqrt{p_{ij}}C^4X_{ij} \)\leq $$
$$\sum_{ij} \(\alpha^2 + \sqrt{p}C^4X_{ij} \)= $$
$$(m-1)^2 \alpha^2 + \sqrt{p}C^4X \leq $$
$$ (m-1)^2\alpha^2 + 4(m-1){\frac{\sqrt\eps}{\alpha}}C^4X$$
This yields
$$X \leq \frac{(m-1)^2\alpha^2}{1-4(m-1)C^4\sqrt{\eps}/\alpha}.$$

\bbox

\subsection{Beckner's inequality}\label{sec:Beckner}
This subsection is an adaptation of results from \cite{Wolff} for our case.

\subsubsection{Introduction to hypercontractivity} 
Let $\Omega$ be an arbitrary finite set endowed with the measure $\mu$. Assume that $\mu$ has at least $2$ atoms with non-zero measure.  Let $\mathfrak{S}$ be a linear subspace of the space of all functions from $\Omega$ to $\R$, that includes the constant functions. 

Define the operator $L$ to be the orthogonal projection to functions with zero mean $L=Id - E_{\mu}$. Define the semigroup $T_t=e^{-tL} (t\geq 0)$. We shall use the explicit formula
$$ T_t=E_{\mu} + e^{-t}L=(1-e^{-t})E_{\mu}+e^{-t}Id$$ 

\dfn
The semigroup $\(T_t\)_{t\geq}$ is $(p,q)$-hypercontractive (for $1<q<p<\infty$) over $\mathfrak{S}$ if there exists $t_0$ such that for all $t \geq t_0$ and $f \in L_q(\mu) \cap \mathfrak{S}$ 
$$\|T_tf\|_{Lp(\mu)} \leq \|f\|_{L_q(\mu)}$$
\edfn
If such a $t_0$ exists and is the least possible, then
$\sigma_{p,q}(\mu,\mathfrak{S})=e^{-t_0}$ is called the
$(p,q)$-\textit{hypercontractive constant} for the measure $\mu$ over
$\mathfrak{S}$.

Multivariate functions are functions from $\Omega^n$ to $\R$. The operators $L$ and $T_t$ have multivariate parallels. Define $$L^i=Id-E_{x_i ~ \mu},\  T^i_t=e^{-tL^i},\ T_t=\prod_{i=1}^n T^i_t$$

From supermultiplicativity of norms, it is easy to observe that for a linear space of functions $  \mathfrak{S} $ as above,
$\sigma_{p,q}(\mu,\mathfrak{S})\leq\sigma_{p,q}(\mu,\mathfrak{S}^n)$. 
However, when $q < p$ one can apply Minkowski's inequality to get the inequality in the opposite direction, and deduce 
$\sigma_{p,q}(\mu,\mathfrak{S})=\sigma_{p,q}(\mu,\mathfrak{S}^n)$.
In our case, $\mathfrak{S}$ is the space of functions  from $\S_m$ to $\R$ whose
Fourier transform is supported on $\rho^0$ and $\rho^1$,

We now wish to consider the degree of multivariate functions. A monomial of degree $d$
is a function $f$ such that there exists a set of
indices $S \subseteq [n]$ of size $d$ such that for $i \in S$, $L^if=f$ and for $i
\notin S$, $L^if=0$.
\dfn \label{dfn:deg} A function of degree $d$ is a function in the
linear span of all monomials of degree $\leq d$ and not in the span of monomials of degree $\leq d-1$.
\edfn

\ni We will use the following corollary of the hypercontractivity for functions of degree $d$:
\clm
If $\sigma_{p,2}(\mu,\mathfrak(S))=e^{-t_0}$, then for $f\in\mathfrak{S}$ of degree $d$, $$ \|f\|_{L_p(\mu)} \leq e^{d \c t_0} \|f\|_{L_2(\mu)}$$
\eclm
\prf
Write $f=T_tT^{-1}_{t}f$. Therefore,
$$\|f\|_{Lp(\mu)} \leq \|T^{-1}_{t}f\|_{L_2(\mu)}$$
Now write $f=f_0+f_1+...+f_d$, where $f_i$ is the projection of $f$ onto the monomials of degree $d$. Then $T^{-1}_{t}f=\sum_{i=1}^d e^{t_o\c i}f_i$. Because the monomials are orthogonal to each other, we have
$$\|T^{-1}_{t}f\|^2_2=\|\sum_{i=1}^d e^{t_o\c i}f_i\|^2_2=\sum_{i=1}^d \|e^{t_o\c i}f_i\|^2_2 \leq e^{2\c t_o\c d} \sum_{i=1}^d \|f_i\|^2_2 = e^{2\c t_o\c d} \|f\|^2_2$$

\eprf
In the next subsubsection, we will calculate $\sigma$ for our needs and arrive at the following lemma:

\cor \label{cor:BeckFinal}
Let $m\geq e^4$. Let $r$ be a function $r:\S_m^n \to \R$, of degree $2$, all of whose Fourier coefficients are supported on tensors of $\rho^0$ and $\rho^1$. Then we have
$$\|r\|_4^4\leq C^8 \|r\|_2^4$$
for $C=\frac 1m ^\frac 12$.
\ecor 

\subsubsection{Calculating $\sigma$} 
In this section, we shall calculate $\sigma=\sigma_{p,q}(\mu,\mathfrak{S_m})$, where $\mathfrak{S}_m$ is the space of functions from $\S_m$ to $\R$, whose Fourier transform is supported on $\rho^0$ and $\rho^1$ (and $\mu$ is the uniform measure). 
Our main result is the following.
\thm \label{thm:SigmaP}
There exists some $m_0$ such that for all $m>m_0$,
$$\sigma_{4,2}(\mu,\mathfrak{S}_m)\geq m^{-\frac 12}$$
\ethm

In the following claim, we present the functions in $\mathfrak{S}_m$ in a form that will help us in the analysis:
\clm \label{clm:Pform}
For $f \in \mathfrak{S}_m$, $f(x)=tr(A\c P(x))$, where $A$ is the matrix $U\(\h{f}(0) \oplus (m-1)\h{f}^t(1)\)U^t$ (See the definition of $U$ in (\ref{eqn:PDec})).
\eclm

\prf
Since $f \in \mathfrak{S}$, we have $f(x)=\h{f}(0)\rho^0(x) + (m-1)tr\(\h{f}^t(1)\rho^1(x)\)$. This equals to
$$
f(x)=\h{f}(0)\rho^0(x) + (m-1)tr\(\h{f}^t(1)\rho^1(x)\) = tr\(\(\h{f}(0) \oplus (m-1)\h{f}^t(1)\) \c \(\rho^0(x) \oplus \rho^1(x)\)\) = $$$$tr\(U\(\h{f}(0) \oplus (m-1)\h{f}^t(1)\)U^{-1} \c U\(\rho^0(x) \oplus \rho^1(x)\)U^{-1}\) =tr(A\c P(x)) 
$$
\eprf
Note that by the definition of $U$, we have $\1 A=(A\1^t)^t=\1\h{f}(0)$.

We define the following moment-like operators on $A$:
$$ M_1(A)=\sum_{i,j}A_{ij}\ \ \ M_2(A)=\sum_{i,j}A^2_{ij}\ \ \ M_3(A)=\sum_{i,j}A^3_{ij}\ \ \ M_4(A)=\sum_{i,j}A^4_{ij}$$$$M_r(A)=\sum_i\(\sum_j A_{ij}^2\)^2\ \ \ M_c(A)=\sum_j\(\sum_i A_{ij}^2\)^2$$$$M_q(A)=tr(AA^tAA^t)=\sum_{ij}\(\sum_k A_{ik}A_{kj}\)^2$$

We can express the $2$'nd and $4$'th norms of functions in $\mathfrak{S}_m$ using these operators, as shown in the following lemmas:
\lem \label{lem:norm2A}
Let $f$ be a function in $\mathfrak{S}_m$ and let $f(x)=tr\(AP(x)\)$ as in claim \ref{clm:Pform}, then
$$|f|_2^2=\E_xf(x)^2=\frac{1}{m-1}\((m-2)M_2(A)+\frac{M_1^2(A)}{m^2}\)$$
\elem
\lem\label{lem:norm4TA}
Let $f$ be a function in $\mathfrak{S}_m$ and let $f(x)=tr\(AP(x)\)$ as in claim \ref{clm:Pform}. For some $t>0$, denote $\sigma=e^{-t}$. Then
$$
\begin{array}{rcl}
|T_tf|_4^4=\E_xf(x)^4=\frac{1}{(m-1)(m-2)(m-3)}
(
& &  \(\frac {m-12}{m^2}+O\(\frac {\sigma^4}{m^3}\)\) M_1^4(A)  \\
&-&  O(\sigma^4) M_1(A)M_3(A)  \\
&-&  O\(\frac{\sigma^4}{m}\) M_1^2(A)M_2(A) \\
&+&  O\(\sigma^4 m\) M_2^2(A)  \\
&+&  \sigma^4\(m^2+m\)M_4(A) \\
&-&  O\(\sigma^4m\) M_r(A) \\
&-&  O\(\sigma^4m\) M_c(A) \\
&+&  O\(\frac{\sigma^4}{m}\) M_q(A))\\
\end{array}$$
\elem
The proofs of these lemmas will be given later.

Next, we shall assume that $f$ is normalized so that that $|f|^2_2=1$ and bound the terms appearing in the expression for $|Tf|^4_4$ under this constraint:

\clm \label{clm:mmntBounds}Let $f$ be a function in $\mathfrak{S}_m$ such that $\|f\|_2^2=1$, then:
\begin{enumerate}
\item $|M_1(A)|\leq O(m) $
\item $M_2(A)\leq m-1$
\item $|M_3{A}|\leq O\(m^{\frac 32}\)$
\item $M_4(A) \leq O(m^2)$
\item $-M_r(A)\leq 0, -M_c(A)\leq 0$
\item $M_q(A)\leq M^2_2(A)\leq O(m^2)$
\end{enumerate}
Also, the leading coefficient in all of these bounds is $1$.
\eclm
Combining all of this information yields Theorem \ref{thm:SigmaP}:

\prf
Assigning the bounds from claim \ref{clm:mmntBounds} in the expression from lemma \ref{lem:norm4TA}, we see that all the term except for three are $o(\frac 1m)$. Addressing the remaining terms we have
$$\|T_tf\|_4^4=\frac{
\frac {m-12}{m^2} M_1^4(A)
+  O\(\sigma^4 m\) M_2^2(A) 
+  \sigma^4m^2M_4(A)}{(m-1)(m-2)(m-3)} +o(\frac 1m)=$$$$ 
\frac{m^2(m-12)
+  O\(\sigma^4 m^3\) 
+  \sigma^4 m^4 }{(m-1)(m-2)(m-3)}
+o(\frac 1m) $$
When choosing $\sigma=m^{-\frac 12}$ we get
$$\|T_tf\|_4^4=\frac{
m^2(m-12)
+  O\(m\)+  m^2 }{(m-1)(m-2)(m-3)}
+o(\frac 1m) $$
This expression is asymptotically smaller than $1$.
\eprf
We shall now present the proofs of the lemmas. We begin with a complicated proof of the following easy lemma.
\lem
Let $f$ be a function in $\mathfrak{S}_m$ and let $f(x)=tr\(AP(x)\)$ as in claim \ref{clm:Pform}, then
$$\E_x f(x)=\frac{M_1(A)}{m}$$
\elem
\prf
$$\E_x f(x)=\E_x tr\(AP(x)\)= tr\(A\E_x P(x)\)$$
Remember that $P$ is a representation and is reducible to a copy of $\rho^0$ and $\rho^1$, so $$E_x P(x)= U\(\E_x \rho^0(x) \oplus \E_x \rho^1(x)\)U^t$$ Because of Schur's orthonormality, we have  $$\E_x \rho^0(x)=<\rho^0,1>=1$$ $$\E_x \rho^1(x)=<\rho^1,1>=0$$ Therefore, $\E_x P(x)= U\(1\oplus 0\)U^t=\frac Jm$. We can now conclude that 
$$\E_x f(x)=tr\(A\frac Jm\)=\frac{M_1(A)}{m}$$
\eprf
We shall use a similar technique to prove lemma \ref{lem:norm2A} regarding the $2$nd norm of $f$.

\ni{\bf Proof of lemma \ref{lem:norm2A}:}We have $$f^2(x)=tr\(AP(x)\) \c tr\(AP(x)\) =tr\(\(A\x A\) \(P(x)\x P(x)\) \)$$and therefore, 
$$\E_x f^2(x) = tr\(\(A\x A\) \E_x\(P(x)\x P(x)\) \)$$
Denote $Q=\E_x [P(x)\x P(x)]$. We need a closed formula for $Q$, which we shall obtain by diagonalizing it. Notice that $P\x P$ is a reducible representation. When taking the expectation $\E_x [P(x)\x P(x)]$, all the copies of $\rho^0$ will become $1$, and the copies of other representations will zero out. The rank of $Q$ is the multiplicity of $\rho^0$ in $P\x P$, which is:
$$<\chi_{\rho^0},\chi_{P\x P}> = \E_x \chi_{P\x P}(x)=\E_x \(\sum_i \sum_j P_{ii}(x)P_{jj}(x)\) = $$$$\sum_i \E_x P_{ii}^2(x) + \sum_{i,j,i\neq j} \E_x P_{ii}(x)P_{jj}(x) = m\c \frac 1m + m(m-1)\c \frac 1{m(m-1)}=2$$
So $Q$'s only nonzero eigenvalue is $1$ with multiplicity $2$. Before computing the eigenvectors of $Q$, notice that $Q$ is symmetric, because $Q^t=\E_x P^t(x)\x P^t(x)=\E_x P(x^{-1})\x P(x^{-1})=Q$. The eigenvectors of $Q$ are of size $m^2$. We can guess two eigenvectors of $Q$. We denote them by $u$ and $v$, and list them by their entries, indexed by $i,j\in[m]$ (Recall that $\1$ is the all ones vector, and $\1_i$ is its $i$th entry, which equals $1$):
$$u_{ij}=\1_i\1_j$$
$$v_{ij}=\delta_{ij}$$ It is easy to see that for every $x$, $u$ and $v$ are eigenvectors of $P(x)\x P(x)$. Therefore, these vectors are also eigenvectors of $Q$. However, these vectors are not orthogonal. We need to find linear combinations of them which are orthonormal. In other words, let $E$ be the $2\times (m^2)$ matrix whose rows are $u$ and $v$. We wish to find a $2 \times 2$ matrix $O$ of coefficients such that $OE$ will be an orthonormal matrix, i.e. $OE\c(OE)^t=OEE^tO^t=I$. In that case, we will have $Q=(OE)^tOE=E^tO^tOE$. Denote $C=EE^t$, then,
$$OCO^t=I\Rightarrow O^{-1}=CO^t \Rightarrow CO^tO=I \Rightarrow O^tO=C^{-1}$$ Since $C=EE^t$ is symmetric and PSD, we may choose $O$ to be $O=C^{-\frac 12}$ and then $O$ is also symmetric.
It is easy to verify that
$$C=EE^t=\(\begin{array}{c c} uu^t&uv^t\\vu^t&vv^t\end{array}\)=\(\begin{array}{c c} m^2&m\\m&m\end{array}\)$$  
and then $$O^tO=OO=C^{-1}=\frac 1{(m-1)m} \(\begin{array}{c c} 1&-1\\-1&m\end{array}\)$$
Now $$tr\(\(A\x A\)Q\)=tr\(\(A\x A\)E^tO^tOE\)=tr\(E\(A\x A\)E^tOO\)=tr\(E\(A\x A\)E^tC^{-1}\)$$
One can verify that 
$$E\(A\x A\)E^t=\(\begin{array}{c c} u\(A\x A\)u^t&u\(A\x A\)v^t\\v\(A\x A\)u^t&v\(A\x A\)v^t\end{array}\)=\(\begin{array}{c c} M_1^2(A)&\frac 1m M_1^2(A)\\\frac 1m M_1^2(A) & M_2(A)\end{array}\)$$
Finally,
$$\E_x f^2(x)=\frac 1{(m-1)m}tr\(\(\begin{array}{c c} M_1^2(A)&\frac 1m M_1^2(A)\\\frac 1m M_1^2(A) & M_2(A)\end{array}\)\c
 \(\begin{array}{c c} 1&-1\\-1&m\end{array}\)\)=$$$$\frac 1{m-1}\(\frac{m-2}{m^2}M_1^2(A) + M_2(A)\)$$\eprf
We shall now use the same technique to compute $\E f^4(x)$:
\lem \label{lem:norm4A}
Let $f$ be a function in $\mathfrak{S}_m$ and let $f(x)=tr\(AP(x)\)$ as in claim \ref{clm:Pform}, then
$$\begin{array}{rcl}
\|f\|_4^4=\E_x f^4(x)=\frac{1}{(m-1)(m-2)(m-3)}
(
& &  \(\frac {m-12}{m^2}+o\(\frac 1{m^2}\)\) M_1^4(A)  \\
&-&  O(1) M_1(A)M_3(A)  \\
&-&  O\(\frac{1}{m}\) M_1^2(A)M_2(A) \\
&+&  O\( m\) M_2^2(A)  \\
&+&  \(m^2+m\)M_4(A) \\
&-&  O\(m\) M_r(A) \\
&-&  O\(m\) M_c(A) \\
&+&  O\(\frac{1}{m}\) M_q(A))\\
\end{array}$$
\elem
\prf
We have $$f^4(x)=tr\(AP(x)\)^4 =tr\(\(A^{\x 4}\) \(P^{\x 4}(x)\) \)$$and therefore, 
$$\E_x f^4(x) = tr\(\(A^{\x 4}\) \E\(P^{\x 4}(x)\) \)$$
Denote $Q=\E_x P^{\x 4}(x)$. We need a closed formula for $Q$, which we shall obtain by diagonalizing it. Notice that $P^{\x 4}$ is a reducible representation. When taking the expectation $\E_x P^{\x 4}(x)$, all the copies of $\rho^0$ will become $1$, and the copies of other representations will zero out. The rank of $Q$ is the multiplicity of $\rho^0$ in $P^{\x 4}$, which is:
$$<\chi_{\rho^0},\chi_{P^{\x 4}}> = \E_x \chi_{P^{\x 4}}(x)=\E_x \(\sum_{i,j,k,l} P_{ii}(x)P_{jj}(x)P_{kk}(x)P_{ll}(x)\) = $$$$\sum_i \E_x P_{ii}^4(x) +4\c\sum_{i,j,i\neq j} \E_x P^3_{ii}(x)P_{jj}(x) +    6\c\sum_{i,j,k,\ distinct} \E_x P^2_{ii}(x)P_{jj}(x)P_{kk}(X) + $$$$3\c\sum_{i,j,i\neq j} \E_x P^2_{ii}(x)P^2_{jj}(x) +  \sum_{i,j,k,l\ distinct} \E_x P_{ii}(x)P_{jj}(x)P_{kk}(x)P_ll(x) = $$$$  \frac mm + 4\c \frac {m(m-1)}{m(m-1)}+6\frac {m(m-1)(m-2)}{m(m-1)(m-2)}+$$$$3\c \frac {m(m-1)}{m(m-1)}+\c \frac {m(m-1)(m-2)(m-3)}{m(m-1)(m-2)(m-3)}=15$$

So $Q$'s only nonzero eigenvalue is $1$ with multiplicity $15$. As before, $Q$ is symmetric. Eigenvectors of $Q$ are of size $m^4$. We now show $15$ eigenvectors of $Q$, divided into $5$ groups. We list them by their entries, indexed by $i,j,k,l \in [m]$:
$$E_1=\{\1_i\1_j\1_k\1_l\}$$
$$E_2=\{\delta_{i,j}\1_k\1_l,\delta_{i,k}\1_j\1_l,\delta_{i,l}\1_j\1_k,\delta_{j,k}\1_i\1_l,\delta_{j,l}\1_i\1_k,\delta_{k,l}\1_i\1_j\}$$
$$E_3=\{\delta_{i,j}\delta_{k,l},\delta_{i,k}\delta_{j,l},\delta_{i,l}\delta_{j,k}\}$$
$$E_4=\{\delta_{i,j,k}\1_l,\delta_{i,j,l}\1_k,\delta_{i,k,l}\1_j,\delta_{j,k,l}\1_i\}$$
$$E_5=\{\delta_{i,j,k,l}\}$$
It is easy to see that for every $x$, These vectors are eigenvectors of $P^{\x 4}(x)$. Therefore, these vectors are also eigenvectors of $Q$. However, these vectors are not orthogonal. As before, let $E$ be the $15\times (m^2)$ matrix whose rows are these vectors. Let $C=EE^t$. Then there exists a symmetric $15 \times 15$ matrix $O$ such that $OE$ is orthonormal, and $OO=C^{-1}$.
In appendix \ref{app:Mat}, we give the expression for $C$ and $E\(A^{\x 4}\)E^t$. We calculated $C^{-1}$ using a mathematical software called \textit{Sage}. The expression is too large to show it here. Using this result, we computed $\E_x f^4(x)=tr\( E\(A^{\x 4}\)E^tC^{-1}\)$ and got the result.

\eprf
We now wish to apply the result of lemma \ref{lem:norm4A} on $T_tf$ to obtain lemma \ref{lem:norm4TA}. We shall use the following two claims:
\clm \label{clm:ATf}
Let $f$ be a function in $\mathfrak{S}_m$ and let $f(x)=tr\(AP(x)\)$ as in claim \ref{clm:Pform}, then $T_tf(x)=tr\(A'P(x)\)$, where $A'=\sigma A+ (1-\sigma)\frac {M_1(A)}{m^2}J$, and $\sigma=e^{-t}$
\eclm
\prf
$$tr\(A'P(x)\)=tr \(\(\sigma A+ (1-\sigma)\frac {M_1(A)}{m^2}J\)P(x)\)=$$$$tr \(\sigma AP(x)\)+ tr\((1-\sigma)\frac {M_1(A)}{m^2}JP(x)\)=\sigma f(x) + (1-\sigma)\E f = T_tf (x)$$
where we have used the fact that $tr(JP(x))=m$ and $M_1(A)=m\c\E f$
\eprf
\clm \label{clm:mmntTA}
For $A$ as in claim \ref{clm:Pform} and $A'=\sigma A+ (1-\sigma)\frac {M_1(A)}{m^2}J$, denote by $\tau=\frac {(1-\sigma)}{m^2}$. We have:
\begin{enumerate}
\item $M_1(A')=M_1(A)$
\item $M_2(A')=\sigma^2 M_2(A)+2\sigma\tau M^2_1(A) + q^2M^2_1(A)$
\item $M_3(A')=\sigma^3 M_3(A)+3\sigma^2\tau M_1(A)M_2(A) +3\sigma\tau^2M_1^3 + \tau^3 M_1^3(A)m^2$
\item $M_4(A')=\sigma^4M_4(A)+4\sigma^3\tau M_1(A)M_3(A)+6\sigma^2\tau^2M_1^2(A)M_2(A)+4\sigma\tau^3M_1^4(A)+\tau^4M_1^4(A)m^2$
\item $M_r(A')=\sigma^4M_r(A) + 4\sigma^3\tau M^2_1(A)M_2(A)\frac 1m + 2\sigma^2\tau^2M_1^2(A)M_2(A)m + 4\sigma^2\tau^2M_1^4(A)\frac 1m + 4\sigma\tau^3M_1^4(A)m + \tau^4M_1^4(A)m^3$
\item $M_c(A')=\sigma^4M_c(A) + 4\sigma^3\tau M^2_1(A)M_2(A)\frac 1m + 2\sigma^2\tau^2M_1^2(A)M_2(A)m + 4\sigma^2\tau^2M_1^4(A)\frac 1m + 4\sigma\tau^3M_1^4(A)m + \tau^4M_1^4(A)m^3$
\item $M_q(A')=\sigma^4M_q(A)+4\sigma^3\tau M_1^4(A)\frac 1{m^2} + 6\sigma^2\tau^2M_1^4(A) + 4\sigma\tau^3M_1^4(A)m^2 + \tau^4M_1^4(A)m^4$
\end{enumerate}
\eclm
\prf
Since all the cases in this lemma are straightforward, we make do with demonstrating the proof of item number $4$, and do not include the tedious, but simple calculations of the other cases. For item $4$ we have:
$$M_4(A')=\sum_{i_1j_1k_1l_1i_2j_2k_2l_2}\delta_{i_1j_1k_1l_1}\(A'\)^{\x 4}_{(i_1j_1k_1l_1)(i_2j_2k_2l_2)}\delta_{i_2j_2k_2l_2}=
$$$$\sum_{i_1j_1k_1l_1i_2j_2k_2l_2}\delta_{i_1j_1k_1l_1}\(\sigma A+ \tau M_1(A)J\)^{\x 4}_{(i_1j_1k_1l_1)(i_2j_2k_2l_2)}\delta_{i_2j_2k_2l_2} 
$$
The above summation can be expanded into six types of summands:
\begin{itemize}
\item $1$ summand of the form
$$\sigma^4\sum_{i_1j_1k_1l_1i_2j_2k_2l_2}\delta_{i_1j_1k_1l_1}\(A\)^{\x 4}_{(i_1j_1k_1l_1)(i_2j_2k_2l_2)}\delta_{i_2j_2k_2l_2}=\sigma^4M_4(A)$$
\item $4$ summands of the form
$$\sigma^3\tau M_1(A)\sum_{i_1j_1k_1l_1i_2j_2k_2l_2}\delta_{i_1j_1k_1l_1}\(A^{\x 3}\x J\)_{(i_1j_1k_1l_1)(i_2j_2k_2l_2)}\delta_{i_2j_2k_2l_2}=\sigma^3\tau M_1(A)M_3(A)$$
\item 
$4$ summands of the form
$$\sigma^2\tau^2 M^2_1(A)\sum_{i_1j_1k_1l_1i_2j_2k_2l_2}\delta_{i_1j_1k_1l_1}\(A^{\x 2}\x J^{\x 2}\)_{(i_1j_1k_1l_1)(i_2j_2k_2l_2)}\delta_{i_2j_2k_2l_2}=$$$$\sigma^2\tau^2 M^2_1(A)M_2(A)$$
\item 
$2$ summands of the form
$$\sigma^2\tau^2 M^2_1(A)\sum_{i_1j_1k_1l_1i_2j_2k_2l_2}\delta_{i_1j_1k_1l_1}\(A\x J\x A\x J\)_{(i_1j_1k_1l_1)(i_2j_2k_2l_2)}\delta_{i_2j_2k_2l_2}=$$$$\sigma^2\tau^2 M^2_1(A)M_2(A)$$
\item 
$4$ summands of the form
$$\sigma\tau^3 M^3_1(A)\sum_{i_1j_1k_1l_1i_2j_2k_2l_2}\delta_{i_1j_1k_1l_1}\(A\x J^{\x 3}\)_{(i_1j_1k_1l_1)(i_2j_2k_2l_2)}\delta_{i_2j_2k_2l_2}=$$$$\sigma\tau^3 M^3_1(A)M_1(A)=\sigma\tau^3 M^4_1(A)$$
\item
$1$ summand of the form
$$\tau^4 M^4_1(A)\sum_{i_1j_1k_1l_1i_2j_2k_2l_2}\delta_{i_1j_1k_1l_1}\(J^{\x 4}\)_{(i_1j_1k_1l_1)(i_2j_2k_2l_2)}\delta_{i_2j_2k_2l_2}=\tau^4 M^4_1(A)m^2$$

\end{itemize} 
\eprf

\ni{\bf Proof of lemma \ref{lem:norm4TA}:} This lemma is simply an application of lemma \ref{lem:norm4A} to $T_tf$ and an assignment of the values of the moments of $A'$ given in claim \ref{clm:mmntTA}. Since there are many terms involved in the assignment, we used \textit{Sage} to get to the actual expression here as well. \qed

\ni{\bf Proof of claim \ref{clm:mmntBounds}:}
\begin{enumerate}
\item $|M_1(A)| \leq O(m) $:$$\E f^2(x)=1=\frac{1}{m-1}\(M_2(A)+(m-1)\frac{M_1^2(A)}{m^2}\)\geq\frac{1}{m-1}\((m-1)\frac{M_1^2(A)}{m^2}\)\Rightarrow $$$$|M_1(A)| \leq O(m)$$
\item $M_2(A)\leq m-1$:
$$\E f^2(x)=1=\frac{1}{m-1}\(M_2(A)+(m-1)\frac{M_1^2(A)}{m^2}\)\geq\frac{1}{m-1}\(M_2(A)\)\Rightarrow $$$$M_2(A)\leq m-1$$
\item $|M_3{A}|\leq O\(m^{\frac 32}\)$: Let $a=\E_{i,j}A_{ij}$, $B=A-aJ$, $b^2=\E B_{i,j}^2$. We take $b$ to be positive. The items above imply that $|a|\leq O(\frac 1m)$ and that $b\leq\sqrt{\frac 1m}$.
$$\sum_{i,j}A_{i,j}^3=\sum_{i,j}(a+B_{i,j})^3=m^2a^3+3aM_2(B)+3a^2M_1(B) +M_3(B)=$$$$m^2a^3+3am^2b^2+M_3(B)$$
Given that the $2$nd moment of $B$ is $m^2b^2$, Jensen's inequality implies that $|M_3(B)|\leq m^3b^3$. Applying the bounds on $a$ and $b$ yields the result.
\item $M_4(A) \leq O(m^2)$: Let $a,B$ and $b$ be as in the previous item. 
$$\sum_{i,j}A_{i,j}^4=\sum_{i,j}(a+B_{i,j})^4=m^2a^4+4aM_3(B)+6a^2M_2(B) + 4a^3M_1(B) +M_4(B)=$$$$m^2a^4+4aM_3(B) + 6a^2m^2b^2+M_4(B)$$
Again, given that the $3$rd moment of $B$ is $m^2b^2$, Jensen's inequality implies that $|M_3(B)|\leq m^3b^3$ and $|M_4(B)|\leq m^4b^4$. Applying the bounds on $a$ and $b$ yields the result.
\item $-M_r(A)\leq 0, -M_c(A)\leq 0$: By definition, these operators are nonnegative.
\item $M_q(A)\leq M^2_2(A)\leq O(m^2)$: Notice that $AA^t$ is a symmetric PSD matrix. Let $\lambda$ be the vector of eigenvalues of $AA^t$ (which are all nonnegative). We have $M_2(A)=tr(AA^t)=\sum_i\lambda_i$ and $M_q(A)=tr\(AA^tAA^t\)=\sum_i\lambda_i^2$. Therefore, $M_q = \sum_i \lambda_i^2 \leq \(\sum_i \lambda_i\)^2= M_2^2$ 
\end{enumerate}
\qed

\section{Strategy Proofness} \label{sec:SP} 
Before discussing further work we wish to briefly discuss the notion
of IR, that arises naturally when considering the spectrum of outcomes
between SCF's and SWF's. The connection between IR and IIA is easy to
understand, despite the fact that neither of these constraints implies
the other. In this section we wish to touch upon the other end of the
spectrum and discuss the connection between IR and strategy proofness.

We present a definition of strategy proofness that is closely related to IR. We will show how our robust impossibility theorem for IR implies a robust impossibility theorem for this definition of strategy proofness. In this sense, our technique provides a single proof for the analogues of both Arrow's theorem (SWF's with independence) and Gibbard-Satterthwaite's theorem (SCF with strategy proofness).

The definition we give here is based on the definition of Dietrich and List \cite{DL} for strategy proofness on judgment aggregation, which is a general framework that includes the setting of Arrow's theorem. In this framework, there is a permissible opinion space $X \subseteq [k]^m$ and we are interested in social aggregators of the form $f:X^n\to X$. A widely studied topic was the characterization of spaces $X$ for which independence implies dictatorship, as in Arrow's theorem and our setting. This characterization can be found in \cite{NP1} and \cite{DH1} for $k=2$ and in \cite{DH} for general $k$. In \cite{DL}, a general definition of strategy proofness was given for this framework was for $k=2$, and was connected to independence. The definition we present here follows from an extension of the definition in \cite{DL} for any $k$, which also deals with the case where the output space of the function is different than the input space. 

Let $H$ be a fixing subgroup of $\S_m$. Let $f$ be $f:\S_m^n \to \S_m/H$. We assume that for any alternative $j \in [m]$, given the ranking of $j$ by a certain voter, we are able to rank the voter's preference over all possible $j$-profiles of the outcome.

Formally, define the set $Q_j$ to be the set of all possible $j$-profiles of a coset of $H$, $Q_j=\{K^{-1}(j) \}_{K \in \S_m/H}$. For every alternative $j \in [m]$ and ranking $r \in [m]$, we assume that there exists some full transitive order relation on $Q_j$. We denote it by $<_{r,j}$. It is natural to assume that for a specific ranking of the alternatives $x$, the top of the order $<_{x^{-1}(j),j}$ will be the $j$-profile of the coset that includes $x$. However, we do not demand that.

A manipulation is a situation where for any alternative $j$, a voter $i\in[n]$ can report a false opinion and get better $j$-profile of the outcome, according to his preference order $<_{x_i^{-1}(j),j}$. Formally, a manipulation exists when the following holds:
$$\exists j,i,x^{-i},x_i,y_i f^{-1}(x^{-i},x_i)(j) <_{x_i^{-1},j} f^{-1}(x^{-i},y_i)(j)$$

The manipulation power of a voter $i$ is the rate of manipulations he can make:
$$M_i(f)=\sum_j \E_{x^{-i},x_i,y_i}\( f^{-1}(x^{-i},x_i)(j) <_{x_i^{-1},j} f^{-1}(x^{-i},y_i)(j)\)$$
The total manipulation power of $f$ is $M(f)=\sum_i M_i(f)$. $f$ will
be called strategy-proof if $M(f)=0$. The rationale behind this definition is that as the designers of a social aggregation mechanism, we do not wish to specify up front for which alternatives $j$ there might be a voter $i$ who wishes to manipulate, and therefore we wish to be immune against all such manipulations.

We shall now connect $IR(f)$ to $M(f)$. Denote $c=max_{K_1,K_2 \in \S_m/H} \left\| \frac{{\bf n}_{K_1}-{\bf n}_{K_2}}{|H|} \right\|_2^2$ (Recall ${\bf n}_K$ is the characteristic vector of a multiset $K$). For example, when $H$ is the group of permutations fixing $\{1\}$, as in SCFs, we have $c=\|(1,0,...0)-(0,\frac 1{m-1},...,\frac 1{m-1})\|_2^2=\frac m{m-1}$.
\clm $cM(f) \geq IR(f)$
\eclm
\prf
Examine all pairs of profiles $(x^{-i},x_i)$ and $(x^{-i},y_i)$, such that $x_i^{-1}(j)=y_i^{-1}(j)$ and $f(x^{-i},x_i)^{-1}(j)\neq f(x^{-i},y_i)^{-1}(j)$. By definition, those pairs of inputs contribute at most $c$ to $IR(f)$.

Denote $r=x_i^{-1}(j)=y_i^{-1}(j)$, $K_x=f(x^{-i},x_i)^{-1}(j)$, $K_y=f(x^{-i},y_i)^{-1}(j)$. Since $K_x\neq K_y$, we must have either $K_x <_{r,j} K_y$ or $K_y <_{r,j} K_x$. Therefore, this pair of inputs contribute $1$ to $M(f)$.
\eprf

As a corollary of this claim and theorem \ref{thm:MainGenRbst}, we have
\cor
\label{Cor:MainSPRbst} Let $H\subseteq \S_m$ be a fixing subgroup of $\S_m$. For $m\geq 3$, an $H$-social aggregator $f$ for which $M(f)\leq\eps$ is $O_H(poly(m) \epsilon)$ close to a function that is either a constant function or dictatorial of the following form: there exists a voter $i$ and a constant permutation $y$ of the rankings such that $f(x)=y\circ x_i$.
\ecor

\section{Further Work}
We have began to apply the techniques of this paper to Arrow's theorem
(with relaxed independence constraint) and to Gibbard-Satterthwaite's
theorem, with partial success. GS seems to offer many more new challenges for this scheme, as the Laplacian is not PSD and not Symmetric.

It is also interesting to see how this technique applies to the
generalized problem of judgment aggregation (see its definition in section \ref{sec:SP}). There are known combinatorial characterization results of functions satisfying independence in this setting, depending on the opinion space $X$, and our technique may be useful in finding robust versions of these theorems, wherever possible.

Another direction is to study how our proof could be modified for groups other than $\S_m$.

As mentioned earlier, our result can be interpreted as a 2 query \textit{dictatorship tester}, for functions $\S_m^n \to \S_m/H$. It is interesting to see whether this has any computational implications.

\bibliographystyle{alpha}
\bibdata{FF}
\bibliography{FF}

\begin{thebibliography}{FKKN11}

\bibitem[Arr63]{Arr}
Kenneth~J. Arrow.
\newblock {\em Social Choice and Individual Values}.
\newblock John Wiley and sons, New York, 1963.

\bibitem[DH10a]{DH1}
Elad Dokow and Ron Holzman.
\newblock Aggregation of binary evaluations.
\newblock {\em Journal of Economic Theory}, 145:495--511, 2010.

\bibitem[DH10b]{DH}
Elad Dokow and Ron Holzman.
\newblock Aggregation of non-binary evaluations.
\newblock {\em Advances in Applied Mathematics}, 45:487--504, 2010.

\bibitem[DL07]{DL}
Franz Dietrich and Christian List.
\newblock Strategy-proof judgement aggregation.
\newblock {\em Economics and Philosophy}, 23:269--300, 2007.

\bibitem[EPF11]{EFP}
Dave Ellis, Haran Pilpel, and Ehud Friedgut.
\newblock Intersecting families of permutations.
\newblock {\em J. Amer. Math. Soc.}, 24:649--682, 2011.

\bibitem[FKKN11]{FKKN}
Ehud Friedgut, Gil Kalai, Nathan Keller, and Noam Nisan.
\newblock A quantitative version of the gibbard-satterthwaite theorem for three
  alternatives.
\newblock {\em SIAM J. Comput.}, 40(3):934--952, 2011.

\bibitem[FKN02]{FKN}
Ehud Friedgut, Gil Kalai, and Assaf Naor.
\newblock Boolean functions whose fourier transform is concentrated on the
  first two levels.
\newblock {\em Adv. in Appl. Math}, 29, 2002.

\bibitem[Gib73]{Gib}
Allan Gibbard.
\newblock Manipulation of voting schemes: a general result.
\newblock {\em Econometrica}, 41(4):587--601, 1973.

\bibitem[IKM10]{IKM}
Marcus Isaksson, Guy Kindler, and Elchanan Mossel.
\newblock The geometry of manipulation: A quantitative proof of the
  gibbard-satterthwaite theorem.
\newblock In {\em FOCS}, pages 319--328. IEEE Computer Society, 2010.

\bibitem[Kal01]{Kal}
Gil Kalai.
\newblock A fourier-theoretic perspective for the condorcet paradox and arrow's
  theorem.
\newblock {\em Advances in Applied Mathematics}, Volume 29, Issue 3:412--426,
  2001.

\bibitem[Kel10a]{Kell2}
Nathan Keller.
\newblock On the probability of a rational outcome for generalized social
  welfare functions on three alternatives.
\newblock {\em Journal of Combinatorial Theory}, Ser. A, 117(4), 2010.

\bibitem[Kel10b]{Kell}
Nathan Keller.
\newblock A tight quantitative version of arrow's impossibility theorem.
\newblock {\em CoRR}, abs/1003.3956, 2010.

\bibitem[Mos09]{M1}
Elchanan Mossel.
\newblock Arrow's impossibility theorem without unanimity.
\newblock {\em CoRR}, abs/0901.4727, 2009.

\bibitem[Mos11]{M}
Elchanan Mossel.
\newblock A quantitative arrow theorem.
\newblock Arxiv 0903.2574, 2011.

\bibitem[MR11]{MR}
Elchanan Mossel and Miklos~Z. Raz.
\newblock A quantitative gibbard-satterthwaite theorem without neutrality.
\newblock {\em arXiv:1110.5888v1}, 2011.

\bibitem[Neh10]{Neh}
Ilan Nehama.
\newblock Approximate judgement aggregation.
\newblock {\em CoRR}, abs/1008.3829, 2010.

\bibitem[NP07]{NP}
Klaus Nehring and Clemens Puppe.
\newblock The structure of strategy-proof social choice part i: General
  characterization and possibility results on median spaces.
\newblock {\em Journal of Economic Theory}, 135:269--305, 2007.

\bibitem[NP10]{NP1}
Klaus Nehring and Clemens Puppe.
\newblock Abstract arrowian aggregation.
\newblock {\em Journal of Economic Theory}, 145:467--494, 2010.

\bibitem[Ren01]{Ren}
Philip~J. Reny.
\newblock Arrow's theorem and the gibbard-satterthwaite theorem: A unified
  approach.
\newblock {\em Economics Letters}, pages 99--105, 2001.

\bibitem[Sat75]{Sat}
Mark~A. Satterthwaite.
\newblock Strategy-proofness and arrow's conditions: Existence and
  correspondence theorems for voting procedures and social welfare functions.
\newblock {\em Journal of Economic Theory}, 10:187--217, April 1975.

\bibitem[Wol07]{Wolff}
Pawel Wolff.
\newblock Hypercontractivity of simple random variables.
\newblock {\em Studia Math.}, 180, iss. 3:219--236, 2007.

\bibitem[Xia08]{Xia}
Lirong Xia.
\newblock A sufficient condition for voting rules to be frequently manipulable.
\newblock In {\em In Proceedings of the Ninth ACM Conference on Electronic
  Commerce (EC}, 2008.

\end{thebibliography}

\appendix
\section{Appendix: matrices for the proof of lemma \ref{lem:norm4A}} \label{app:Mat}
$$C=EE^t=\(\begin{array}{c |c c c c c c |c c c |c c c c |c } 
m^4& m^3&m^3&m^3&m^3&m^3&m^3& m^2&m^2&m^2& m^2&m^2&m^2&m^2& m^1\\
\hline
m^3& m^3&m^2&m^2&m^2&m^2&m^2& m^2&m^1&m^1& m^2&m^2&m^1&m^1& m^1\\
m^3& m^2&m^3&m^2&m^2&m^2&m^2& m^1&m^2&m^1& m^2&m^1&m^2&m^1& m^1\\
m^3& m^2&m^2&m^3&m^2&m^2&m^2& m^1&m^1&m^2& m^1&m^2&m^2&m^1& m^1\\
m^3& m^2&m^2&m^2&m^3&m^2&m^2& m^1&m^1&m^2& m^2&m^1&m^1&m^2& m^1\\
m^3& m^2&m^2&m^2&m^2&m^3&m^2& m^1&m^2&m^1& m^1&m^2&m^1&m^2& m^1\\
m^3& m^2&m^2&m^2&m^2&m^2&m^3& m^2&m^1&m^1& m^1&m^1&m^2&m^2& m^1\\
\hline
m^2& m^2&m^1&m^1&m^1&m^1&m^2& m^2&m^1&m^1& m^1&m^1&m^1&m^1& m^1\\
m^2& m^1&m^2&m^1&m^1&m^2&m^1& m^1&m^2&m^1& m^1&m^1&m^1&m^1& m^1\\
m^2& m^1&m^1&m^2&m^2&m^1&m^1& m^1&m^1&m^2& m^1&m^1&m^1&m^1& m^1\\
\hline
m^2& m^2&m^2&m^1&m^2&m^1&m^1& m^1&m^1&m^1& m^2&m^1&m^1&m^1& m^1\\
m^2& m^2&m^1&m^2&m^1&m^2&m^1& m^1&m^1&m^1& m^1&m^2&m^1&m^1& m^1\\
m^2& m^1&m^2&m^2&m^1&m^1&m^2& m^1&m^1&m^1& m^1&m^1&m^2&m^1& m^1\\
m^2& m^1&m^1&m^1&m^2&m^2&m^2& m^1&m^1&m^1& m^1&m^1&m^1&m^2& m^1\\
\hline
m^1& m^1&m^1&m^1&m^1&m^1&m^1& m^1&m^1&m^1& m^1&m^1&m^1&m^1& m^1
\end{array}\)$$  
We compute $C^{-1}$ using an algebraic software called \textit{Sage}. The expressions are too large to write here. Since all the entries of $C$ are monomials in $m$, the determinant of $C$ is a polynomial of a bounded degree in $m$, and therefore $C$ is regular for all but a finite number of values of $m$. Actually, using \textit{Sage}, we find the determinant of $C$ to be $m^{15}(m-1)^{14}(m-2)^{7}(m-3)$, so $C$ is singular only for $m=0,1,2,3$.

For the expression for $E\(A^{\x 4}\)E$ we abuse the notation and use the sets $E_i$ as matrices whose rows are the elements of $E_i$. We have:
$$E\(A^{\x 4}\)E = \(\begin{array}{c c c c}
E_1\(A^{\x 4}\)E_1^t& E_1\(A^{\x 4}\)E_2^t& \ldots&  E_1\(A^{\x 4}\)E_5^t\\ 
\vdots&\vdots&\ddots&\vdots\\
E_5\(A^{\x 4}\)E_1^t& E_5\(A^{\x 4}\)E_2^t& \ldots&  E_5\(A^{\x 4}\)E_5^t\\ 
\end{array}\)$$
Where:

$$
E_1\(A^{\x 4}\)E_1^t = 
\(\begin{array}{c }
M_1^4 \end{array}\right)
$$$$
E_1\(A^{\x 4}\)E_2^t = 
\frac{1}{m}\(\begin{array}{c c c c c c}
M_1^4&M_1^4&M_1^4&M_1^4&M_1^4&M_1^4 \end{array}\right)
$$$$
E_1\(A^{\x 4}\)E_3^t = 
m\(\begin{array}{c c c}
M_1^4&M_1^4&M_1^4 \end{array}\right)
$$$$
E_1\(A^{\x 4}\)E_4^t = 
m\(\begin{array}{c c c c }
M_1^4&M_1^4&M_1^4&M_1^4 \end{array}\right)
$$$$
E_1\(A^{\x 4}\)E_5^t = 
\frac{1}{m^2}\(\begin{array}{c }
M_1^4\\
\end{array}\right)
$$$$
E_2\(A^{\x 4}\)E_2^t = 
\frac{1}{m^2}\(\begin{array}{c c c c c c}
M^2_1M_2&M_1^4m&M_1^4m&M_1^4m&M_1^4m&M_1^4m\\
M_1^4m&M^2_1M_2&M_1^4m&M_1^4m&M_1^4m&M_1^4m\\
M_1^4m&M_1^4m&M^2_1M_2&M_1^4m&M_1^4m&M_1^4m\\
M_1^4m&M_1^4m&M_1^4m&M^2_1M_2&M_1^4m&M_1^4m\\
M_1^4m&M_1^4m&M_1^4m&M_1^4m&M^2_1M_2&M_1^4m\\
M_1^4m&M_1^4m&M_1^4m&M_1^4m&M_1^4m&M^2_1M_2
\end{array}\right)
$$$$
E_2\(A^{\x 4}\)E_3^t = \frac{1}{m}
\(\begin{array}{c c c }
M^2_1M_2&M_1^4\frac{1}{m^2}&M_1^4\frac{1}{m^2}\\
M_1^4\frac{1}{m^2}&M^2_1M_2&M_1^4\frac{1}{m^2}\\
M_1^4\frac{1}{m^2}&M_1^4\frac{1}{m^2}&M^2_1M_2\\
M_1^4\frac{1}{m^2}&M_1^4\frac{1}{m^2}&M^2_1M_2\\
M_1^4\frac{1}{m^2}&M^2_1M_2&M_1^4\frac{1}{m^2}\\
M^2_1M_2&M_1^4\frac{1}{m^2}&M_1^4\frac{1}{m^2}
\end{array}\right)
$$$$
E_2\(A^{\x 4}\)E_4^t = \frac{1}{m}
\(\begin{array}{c c c c }
M^2_1M_2&M^2_1M_2&M_1^4\frac{1}{m^2}&M_1^4\frac{1}{m^2}\\
M^2_1M_2&M_1^4\frac{1}{m^2}&M^2_1M_2&M_1^4\frac{1}{m^2}\\
M_1^4\frac{1}{m^2}&M^2_1M_2&M^2_1M_2&M_1^4\frac{1}{m^2}\\
M^2_1M_2&M_1^4\frac{1}{m^2}&M_1^4\frac{1}{m^2}&M^2_1M_2\\
M_1^4\frac{1}{m^2}&M^2_1M_2&M_1^4\frac{1}{m^2}&M^2_1M_2\\
M_1^4\frac{1}{m^2}&M_1^4\frac{1}{m^2}&M^2_1M_2&M^2_1M_2\\
\end{array}\right)
$$$$
E_2\(A^{\x 4}\)E_5^t = \frac{1}{m^2}
\(\begin{array}{c }
M^2_1M_2\\x^2_1M_2\\x^2_1M_2\\x^2_1M_2\\
\end{array}\right)
$$$$
E_3\(A^{\x 4}\)E_3^t = 
\(\begin{array}{c c c}
M^2_2&M_q&M_q\\
M_q&M^2_2&M_q\\
M_q&M_q&M^2_2
\end{array}\right)
$$$$
E_3\(A^{\x 4}\)E_4^t = 
\frac{1}{m^2}\(\begin{array}{c c c c }
M^2_1M_2&M^2_1M_2&M^2_1M_2&M^2_1M_2\\
M^2_1M_2&M^2_1M_2&M^2_1M_2&M^2_1M_2\\
M^2_1M_2&M^2_1M_2&M^2_1M_2&M^2_1M_2
\end{array}\right)
$$$$
E_3\(A^{\x 4}\)E_5^t = 
\(\begin{array}{c}
M_c\\M_c\\M_c\end{array}\right)
$$$$
E_4\(A^{\x 4}\)E_4^t = 
\frac{1}{m^2}\(\begin{array}{c c c c }
M_1M_3\frac{1}{m}&M^2_1M_2&M^2_1M_2&M^2_1M_2\\
M^2_1M_2&M_1M_3\frac{1}{m}&M^2_1M_2&M^2_1M_2\\
M^2_1M_2&M^2_1M_2&M_1M_3\frac{1}{m}&M^2_1M_2\\
M^2_1M_2&M^2_1M_2&M^2_1M_2&M_1M_3\frac{1}{m}
\end{array}\right)
$$$$
E_4\(A^{\x 4}\)E_5^t = 
\frac{1}{m}\(\begin{array}{c }
M_1M_3\\M_1M_3\\M_1M_3\\M_1M_3
\end{array}\right)
$$$$
E_5\(A^{\x 4}\)E_3^t = 
\(\begin{array}{c c c}
M_c&M_c&M_c
\end{array}\right)
$$$$
E_5\(A^{\x 4}\)E_5^t = 
\(\begin{array}{c }
M_4
\end{array}\right)
$$
For the blocks that weren't explicitly mentioned, $E_i\(A^{\x 4}\)E_j^t=\(E_i\(A^{\x 4}\)E_j^t\)^t$.

\ni We add here the protocol of the Sage code we used, for completeness:

\ttfamily
R.<m,a1,a13,a112,a22,a4,r22,c22,q,sig>=QQ[]

C=matrix(R,15,15,

[

m\^{}4, m\^{}3,m\^{}3,m\^{}3,m\^{}3,m\^{}3,m\^{}3, m\^{}2,m\^{}2,m\^{}2, m\^{}2,m\^{}2,m\^{}2,m\^{}2, m\^{}1,

m\^{}3, m\^{}3,m\^{}2,m\^{}2,m\^{}2,m\^{}2,m\^{}2, m\^{}2,m\^{}1,m\^{}1, m\^{}2,m\^{}2,m\^{}1,m\^{}1, m\^{}1,

m\^{}3, m\^{}2,m\^{}3,m\^{}2,m\^{}2,m\^{}2,m\^{}2, m\^{}1,m\^{}2,m\^{}1, m\^{}2,m\^{}1,m\^{}2,m\^{}1, m\^{}1,

m\^{}3, m\^{}2,m\^{}2,m\^{}3,m\^{}2,m\^{}2,m\^{}2, m\^{}1,m\^{}1,m\^{}2, m\^{}1,m\^{}2,m\^{}2,m\^{}1, m\^{}1,

m\^{}3, m\^{}2,m\^{}2,m\^{}2,m\^{}3,m\^{}2,m\^{}2, m\^{}1,m\^{}1,m\^{}2, m\^{}2,m\^{}1,m\^{}1,m\^{}2, m\^{}1,

m\^{}3, m\^{}2,m\^{}2,m\^{}2,m\^{}2,m\^{}3,m\^{}2, m\^{}1,m\^{}2,m\^{}1, m\^{}1,m\^{}2,m\^{}1,m\^{}2, m\^{}1,

m\^{}3, m\^{}2,m\^{}2,m\^{}2,m\^{}2,m\^{}2,m\^{}3, m\^{}2,m\^{}1,m\^{}1, m\^{}1,m\^{}1,m\^{}2,m\^{}2, m\^{}1,

m\^{}2, m\^{}2,m\^{}1,m\^{}1,m\^{}1,m\^{}1,m\^{}2, m\^{}2,m\^{}1,m\^{}1, m\^{}1,m\^{}1,m\^{}1,m\^{}1, m\^{}1,

m\^{}2, m\^{}1,m\^{}2,m\^{}1,m\^{}1,m\^{}2,m\^{}1, m\^{}1,m\^{}2,m\^{}1, m\^{}1,m\^{}1,m\^{}1,m\^{}1, m\^{}1,

m\^{}2, m\^{}1,m\^{}1,m\^{}2,m\^{}2,m\^{}1,m\^{}1, m\^{}1,m\^{}1,m\^{}2, m\^{}1,m\^{}1,m\^{}1,m\^{}1, m\^{}1,

m\^{}2, m\^{}2,m\^{}2,m\^{}1,m\^{}2,m\^{}1,m\^{}1, m\^{}1,m\^{}1,m\^{}1, m\^{}2,m\^{}1,m\^{}1,m\^{}1, m\^{}1,

m\^{}2, m\^{}2,m\^{}1,m\^{}2,m\^{}1,m\^{}2,m\^{}1, m\^{}1,m\^{}1,m\^{}1, m\^{}1,m\^{}2,m\^{}1,m\^{}1, m\^{}1,

m\^{}2, m\^{}1,m\^{}2,m\^{}2,m\^{}1,m\^{}1,m\^{}2, m\^{}1,m\^{}1,m\^{}1, m\^{}1,m\^{}1,m\^{}2,m\^{}1, m\^{}1,

m\^{}2, m\^{}1,m\^{}1,m\^{}1,m\^{}2,m\^{}2,m\^{}2, m\^{}1,m\^{}1,m\^{}1, m\^{}1,m\^{}1,m\^{}1,m\^{}2, m\^{}1,

m\^{}1, m\^{}1,m\^{}1,m\^{}1,m\^{}1,m\^{}1,m\^{}1, m\^{}1,m\^{}1,m\^{}1, m\^{}1,m\^{}1,m\^{}1,m\^{}1, m\^{}1

])

t=(1-sig)/m\^{}2

\# the variables a1,a13,a112,a22,a4,r22,c22,q represent, repectively $$M_1^4(A), M_1(A)M_3(A), M_1^2(A)M_2(A), M_2^2(A), M_4(A), M_r(A), M_c(A), M_q(A)$$ 

\# the variables a,b,c,d,e,f,g,h represent, repectively
$$M_1^4(A'), M_1(A')M_3(A'), M_1^2(A')M_2(A'), M_2^2(A'), M_4(A'), M_r(A'), M_c(A'), M_q(A')$$ 

a=a1

b=sig\^{}3*a13 + 3*sig*(1-sig)\^{}2*a1/m\^{}4 + 3*sig\^{}2*(1-sig)*a112/m\^{}2 + 
t\^{}3*a1*m\^{}2

c=sig\^{}2*a112 + (2*sig*t+t\^{}2)*a1

d=sig\^{}4*a22 + 2*(2*sig*t+t\^{}2)*sig\^{}2*a112 + (2*sig*t+t\^{}2)\^{}2*a1 

e=sig\^{}4*a4 + 4*sig*t\^{}3*a1 + 6*sig\^{}2*t\^{}2*a112 + 4*sig\^{}3*t*a13 + t\^{}4*a1*m\^{}2

f=sig\^{}4*r22 + 4*sig\^{}3*t*a112/m + sig\^{}2*t\^{}2*(2*a112*m + 4*a1/m) + 4*sig*t\^{}3*a1*m + t\^{}4*a1*m\^{}3

g=sig\^{}4*c22 + 4*sig\^{}3*t*a112/m + sig\^{}2*t\^{}2*(2*a112*m + 4*a1/m) + 4*sig*t\^{}3*a1*m + t\^{}4*a1*m\^{}3

h=sig\^{}4*q + 4*sig\^{}3*t*a1/m\^{}2 + 6*sig\^{}2*t\^{}2*a1 + 4*sig*t\^{}3*a1*m\^{}2 + t\^{}4*a1*m\^{}4

\#Sage can not work with matrices with rational entries. Therefore, we multiply all the variables by m\^{}10 and the matrix E by m\^{}3.

a=a*m\^{}10; b=b*m\^{}10; c=c*m\^{}10; d=d*m\^{}10; e=e*m\^{}10; f=f*m\^{}10; g=g*m\^{}10; h=h*m\^{}10;

E=matrix(R,15,15,

[

a*m\^{}3, a*m\^{}2,a*m\^{}2,a*m\^{}2,a*m\^{}2,a*m\^{}2,a*m\^{}2, a*m\^{}1,a*m\^{}1,a*m\^{}1, a*m\^{}1,a*m\^{}1,a*m\^{}1, a*m\^{}1, a,

a*m\^{}2, c,a*m,a*m,a*m,a*m,a*m, c*m\^{}2,a,a, c*m\^{}2,c*m\^{}2,a,a, c*m,

a*m\^{}2, a*m,c,a*m,a*m,a*m,a*m, a,c*m\^{}2,a, c*m\^{}2,a,c*m\^{}2,a, c*m,

a*m\^{}2, a*m,a*m,c,a*m,a*m,a*m, a,a,c*m\^{}2, a,c*m\^{}2,c*m\^{}2,a, c*m,

a*m\^{}2, a*m,a*m,a*m,c,a*m,a*m, a,a,c*m\^{}2, c*m\^{}2,a,a,c*m\^{}2, c*m,

a*m\^{}2, a*m,a*m,a*m,a*m,c,a*m, a,c*m\^{}2,a, a,c*m\^{}2,a,c*m\^{}2, c*m,

a*m\^{}2, a*m,a*m,a*m,a*m,a*m,c, c*m\^{}2,a,a, a,a,c*m\^{}2,c*m\^{}2, c*m,

a*m\^{}1, c*m\^{}2,a,a,a,a,c*m\^{}2, d*m\^{}3,h*m\^{}3,h*m\^{}3, c*m,c*m,c*m,c*m, f*m\^{}3,

a*m\^{}1, a,c*m\^{}2,a,a,c*m\^{}2,a, h*m\^{}3,d*m\^{}3,h*m\^{}3, c*m,c*m,c*m,c*m, f*m\^{}3,

a*m\^{}1, a,a,c*m\^{}2,c*m\^{}2,a,a, h*m\^{}3,h*m\^{}3,d*m\^{}3, c*m,c*m,c*m,c*m, f*m\^{}3,

a*m\^{}1, c*m\^{}2,c*m\^{}2,a,c*m\^{}2,a,a, c*m,c*m,c*m, b,c*m,c*m,c*m, b*m\^{}2,

a*m\^{}1, c*m\^{}2,a,c*m\^{}2,a,c*m\^{}2,a, c*m,c*m,c*m, c*m,b,c*m,c*m, b*m\^{}2,

a*m\^{}1, a,c*m\^{}2,c*m\^{}2,a,a,c*m\^{}2, c*m,c*m,c*m, c*m,c*m,b,c*m, b*m\^{}2,

a*m\^{}1, a,a,a,c*m\^{}2,c*m\^{}2,c*m\^{}2, c*m,c*m,c*m, c*m,c*m,c*m,b, b*m\^{}2,

a,     c*m,c*m,c*m,c*m,c*m,c*m, g*m\^{}3,g*m\^{}3,g*m\^{}3, b*m\^{}2,b*m\^{}2,b*m\^{}2,b*m\^{}2, e*m\^{}3

])

l=(C.inverse()*E).trace()/m\^{}13

l(m,1,0,0,0,0,0,0,0,sig)

l(m,0,1,0,0,0,0,0,0,sig)

l(m,0,0,1,0,0,0,0,0,sig)

l(m,0,0,0,1,0,0,0,0,sig)

l(m,0,0,0,0,1,0,0,0,sig)

l(m,0,0,0,0,0,1,0,0,sig)

l(m,0,0,0,0,0,0,1,0,sig)

l(m,0,0,0,0,0,0,0,1,sig)

\end{document}